\documentclass{amsart}%
\usepackage{graphicx}
\usepackage{amscd}
\usepackage{amsmath}
\usepackage{amsfonts}
\usepackage{amssymb}%
\newtheorem{theorem}{Theorem}[section]

\newtheorem{corollary}[theorem]{Corollary}

\newtheorem{definition}[theorem]{Definition}

\newtheorem{lemma}[theorem]{Lemma}

\newtheorem{proposition}[theorem]{Proposition}
\newtheorem{remark}[theorem]{Remark}

\numberwithin{equation}{section}

\begin{document}
\title[Cobordisms of maps]{Cobordisms of maps without prescribed singularities}
\author{YOSHIFUMI ANDO}
\address{Department of Mathematical Sciences, Faculty of Science, Yamaguchi University,
Yamaguchi 753-8512, Japan}
\email{andoy@yamaguchi-u.ac.jp}
\thanks{2000 \textit{Mathematics Subject Classification.} 57R45, 57R90, 58A20}
\date{}
\dedicatory{ }
\begin{abstract}
Let $N$ and $P$ be smooth closed manifolds of dimensions $n$ and $p$
respectively. Given a Thom-Boardman symbol $I$, a smooth map $f:N\rightarrow
P$ is called an $\Omega^{I}$-regular map if and only if the Thom-Boardman
symbol of each singular point of $f$ is not greater than $I$ in the
lexicographic order. We will represent the group of all cobordism classes of
$\Omega^{I}$-regular maps of $n$-dimensional closed manifolds into $P$ in
terms of certain stable homotopy groups. As an application we will study the
relationship among the stable homotopy groups of spheres, the above cobordism
group and higher singularities.

\end{abstract}
\maketitle

\section{Introduction}

Let $N$ and $P$ be smooth ($C^{\infty}$) manifolds of dimensions $n$ and $p$
respectively. Let $k\gg n,$ $p$. Let $J^{k}(N,P)$ denote the $k$-jet bundle of
the manifolds $N$ and $P$ with the canonical projection $\pi_{N}^{k}\times
\pi_{P}^{k}$ onto $N\times P$, whose fiber is denoted by $J^{k}(n,p)$. Here,
$\pi_{N}^{k}$ and $\pi_{P}^{k}$ map a $k$-jet\ to its source and target
respectively. Let $I=(i_{1},i_{2},\cdots,i_{k})$ be a Thom-Boardman symbol
(simply symbol) where $i_{1},i_{2},\cdots,i_{k}$ are a finite number of
integers with $i_{1}\geq i_{2}\geq\cdots\geq i_{k}\geq0$. In \cite{Board}
there have been defined what is called the Boardman manifold $\Sigma^{I}(N,P)$
in $J^{k}(N,P)$. A smooth map germ $f:(N,x)\rightarrow(P,y)$ has $x$ as a
singularity of the symbol $I$ if and only if $j_{x}^{k}f\in\Sigma^{I}(N,P)$.
Let $\Omega^{I}(N,P)$ denote the open subset of $J^{k}(N,P)$ which consists of
all Boardman manifolds $\Sigma^{I^{\prime}}(N,P)$ with symbols $I^{\prime}$ of
length $k$ and $I^{\prime}\leq I$. It is known that $\Omega^{I}(N,P)$ is an
open subbundle of $J^{k}(N,P)$ over $N\times P$, whose fiber is denoted by
$\Omega^{I}(n,p)$. A smooth map $f:N\rightarrow P$ is called an $\Omega^{I}%
$-\textit{regular map} if and only if $j^{k}f(N)\subset\Omega^{I}(N,P).$

Let $J$ be another Thom-Boardman symbol with $I\leq J$. In this paper we will
represent the set of all cobordism classes as $\Omega^{J}$-regular maps of
$\Omega^{I}$-regular maps of $n$-dimensional closed manifolds into $P$ in
terms of certain stable homotopy groups.

Let $P$ be a connected (resp. an oriented) smooth manifold of dimension $p$.
We define the notion of (resp. oriented) $\Omega^{J}$-cobordisms of
$\Omega^{I}$-regular maps. Let $f_{i}:N_{i}\rightarrow P$ ($i=0,1$) be two
$\Omega^{I}$-regular maps, where $N_{i}$ are closed (resp. oriented) smooth
$n$-dimensional manifolds. We say that they are\emph{\ }(resp. \emph{oriented}%
)\emph{\ }$\Omega^{J}$\emph{-cobordant} when there exists an $\Omega^{J}%
$-regular map, say $\Omega^{J}$-cobordism $F:(W,\partial W)\rightarrow
(P\times\lbrack0,1],P\times0\cup P\times1)$ such that, for a sufficiently
small positive number $\varepsilon$,

(i) $W$ is a (resp. an oriented) smooth manifold of dimension $n+1$ with
$\partial W=N_{0}\cup(-N_{1})$ and the collar of $\partial W$ is identified
with $N_{0}\times\lbrack0,\varepsilon]\cup N_{1}\times\lbrack1-\varepsilon,1]$,

(ii) $F|N_{0}\times\lbrack0,\varepsilon]=f_{0}\times id_{[0,\varepsilon]}$ and
$F|N_{1}\times\lbrack1-\varepsilon,1]=f_{1}\times id_{[1-\varepsilon,1]}.$

\noindent Let $\mathfrak{NCob}_{n,P}^{(\Omega_{n}^{I},\Omega_{n+1}^{J})}$
(resp. $\mathfrak{OCob}_{n,P}^{(\Omega_{n}^{I},\Omega_{n+1}^{J})}$) denote the
set of all (resp. oriented) $\Omega^{J}$-cobordism classes of $\Omega^{I}%
$-regular maps into $P$. We will provide them the structures of modules in
Section 4.

We need some notion to represent them by using stable homotopy groups. Let
$E\rightarrow X$ and $F\rightarrow Y\,$be smooth vector bundles of dimensions
$n$ and $p$ over smooth manifolds, and let $\pi_{X}$ and $\pi_{Y}$ be the
projections of $X\times Y$ onto $X$ and $Y$ respectively.\ Define the vector
bundle $J^{k}(E,F)$ over $X\times Y$ by%
\begin{equation}
J^{k}(E,F)=\bigoplus_{i=1}^{k}\text{\textrm{Hom}}(S^{i}(\pi_{X}^{\ast}%
(E)),\pi_{Y}^{\ast}(F)).
\end{equation}
Here, $S^{i}(E)$ is the vector bundle $\cup_{x\in X}S^{i}(E_{x})$ over $X$,
where $S^{i}(E_{x})$ denotes the $i$-fold symmetric product of $E_{x}$. The
canonical fiber $\bigoplus_{i=1}^{k}\mathrm{Hom}(S^{i}(\mathbb{R}%
^{n}),\mathbb{R}^{p})$ is canonically identified with $J^{k}(n,p)$. If we
provide $N$ and $P$ with Riemannian metrics, then $J^{k}(TN,TP)$ is identified
with $J^{k}(N,P)$\ over $N\times P$ (see Section 2). Let $\Omega^{I}(E,F)$
denote the open subbundle of $J^{k}(E,F)$ associated to $\Omega^{I}(n,p)$.

Let $G_{m}$ refer to the grassmann manifold $G_{m,\ell}$\ (resp. oriented
grassmann manifold $\widetilde{G}_{m,\ell}$) of all (resp. oriented)
$m$-subspaces of $\mathbb{R}^{m+\ell}$.\ Let $\gamma_{G_{m}}^{m}$ and
$\widehat{\gamma}_{G_{m}}^{\ell}$\ denote the canonical vector bundles of
dimensions $m$ and $\ell$\ over the space $G_{m}$\ respectively such that
$\gamma_{G_{m}}^{m}\oplus\widehat{\gamma}_{G_{m}}^{\ell}$ is the trivial
bundle $\varepsilon_{G_{m}}^{m+\ell}$. Let $T(\widehat{\gamma}_{G_{m}}^{\ell
})$ denote the Thom spaces of $\widehat{\gamma}_{G_{m}}^{\ell}$. Let
$i^{G}:G_{n}\rightarrow G_{n+1}$ denote the injection mapping an $n$-plane $a$
to the $(n+1)$-plane generated by $a$ and the $(n+\ell+1)$-th unit vector
$\mathbf{e}_{n+\ell+1}$ in $\mathbb{R}^{n+\ell+1}$. Let $i_{+1}^{(\Omega
_{n}^{I},\Omega_{n+1}^{J})}:\Omega^{I}(n,p)\rightarrow\Omega^{J}(n+1,p+1)$
denote the map which sends a $k$-jet $z=j_{0}^{k}f$ to $j_{0}^{k}(f\times
id_{\mathbb{R}})$, where $f:\mathbb{R}^{n}\rightarrow\mathbb{R}^{p}$ is a map
germ and $id_{\mathbb{R}}$ is the identity of $\mathbb{R}$. Set
$\boldsymbol{\Omega}_{n}^{I}=\Omega^{I}(\gamma_{G_{n}}^{n},TP)$,
$\boldsymbol{\Omega}_{n+1}^{J}=\Omega^{J}(\gamma_{G_{n+1}}^{n+1}%
,TP\oplus\varepsilon_{P}^{1})$, $\widehat{\gamma}_{\boldsymbol{\Omega}_{n}%
^{I}}^{\ell}=(\pi_{G_{n}}^{k})^{\ast}(\widehat{\gamma}_{G_{n}}^{\ell
})|_{\boldsymbol{\Omega}_{n}^{I}}$ and $\widehat{\gamma}_{\boldsymbol{\Omega
}_{n+1}^{J}}^{\ell}=(\pi_{G_{n+1}}^{k})^{\ast}(\widehat{\gamma}_{G_{n+1}%
}^{\ell})|_{\boldsymbol{\Omega}_{n+1}^{J}}$, where $\varepsilon_{P}%
^{1}=P\times\mathbb{R}$. There exists a fiberwise map $\mathbf{j}^{(\Omega
_{n}^{I},\Omega_{n+1}^{J})}:\boldsymbol{\Omega}_{n}^{I}\rightarrow
\boldsymbol{\Omega}_{n+1}^{J}$ associated to $i_{+1}^{(\Omega_{n}^{I}%
,\Omega_{n+1}^{J})}$ covering $i^{G}\times id_{P}$. Then $\mathbf{j}%
^{(\Omega_{n}^{I},\Omega_{n+1}^{J})}$ induces the bundle map
\[
\mathbf{b(}\widehat{\gamma})^{(\Omega_{n}^{I},\Omega_{n+1}^{J})}%
:\widehat{\gamma}^{\ell}{}_{\boldsymbol{\Omega}_{n}^{I}}\longrightarrow
\widehat{\gamma}_{\boldsymbol{\Omega}_{n+1}^{J}}^{\ell}%
\]
covering $\mathbf{j}^{(\Omega_{n}^{I},\Omega_{n+1}^{J})}$ and the Thom map
$T(\mathbf{b(}\widehat{\gamma})^{(\Omega_{n}^{I},\Omega_{n+1}^{J})})$. We
denote the image of%
\begin{equation}
T(\mathbf{b(}\widehat{\gamma})^{(\Omega_{n}^{I},\Omega_{n+1}^{J})})_{\ast}%
:\pi_{n+\ell}\left(  \widehat{\gamma}_{\boldsymbol{\Omega}_{n}^{I}}^{\ell
}\right)  \longrightarrow\pi_{n+\ell}\left(  \widehat{\gamma}%
_{\boldsymbol{\Omega}_{n+1}^{J}}^{\ell}\right)
\end{equation}
by $\mathrm{\operatorname{Im}}^{\mathfrak{N}}\left(  T(\mathbf{b(}%
\widehat{\gamma})^{(\Omega_{n}^{I},\Omega_{n+1}^{J})})\right)  $ (resp.
$\mathrm{\operatorname{Im}}^{\mathfrak{O}}\left(  T(\mathbf{b(}\widehat
{\gamma})^{(\Omega_{n}^{I},\Omega_{n+1}^{J})})\right)  $) in the nonoriented
(resp.\ the oriented) case. We are ready to state the main result of this
paper. The following theorem will be proved by applying the homotopy
principles in \cite{duPMCS}, \cite{duPRS} and \cite{hPrinBo}.

\begin{theorem}
Let $n$ and $p$ be natural numbers with $p>1$ and $\ell\gg n$,$p$. Let $P$ be
a connected $p$-dimensional manifold and be oriented in the oriented case. Let
$I$ and $J$ be Thom-Boardman symbols with $I\leq J$ such that if $n\geq p$,
then $I\geq(n-p+1,0)$. Then there exist the isomorphisms%
\begin{align*}
\omega_{\mathfrak{N}}^{(\Omega_{n}^{I},\Omega_{n+1}^{J})}  &
:\mathfrak{{NCob}}_{n,P}^{(\Omega_{n}^{I},\Omega_{n+1}^{J})}\longrightarrow
\mathrm{\operatorname{Im}}\mathfrak{^{\mathfrak{N}}}(T(\mathbf{b(}%
\widehat{\gamma})^{(\Omega_{n}^{I},\Omega_{n+1}^{J})}))\mathfrak{,}\\
\omega_{\mathfrak{O}}^{(\Omega_{n}^{I},\Omega_{n+1}^{J})}  &
:\mathfrak{{OCob}}_{n,P}^{(\Omega_{n}^{I},\Omega_{n+1}^{J})}\longrightarrow
\mathrm{\operatorname{Im}}\mathfrak{^{\mathfrak{O}}}(T(\mathbf{b(}%
\widehat{\gamma})^{(\Omega_{n}^{I},\Omega_{n+1}^{J})})).
\end{align*}

\end{theorem}

As for the image of $T(\mathbf{b(}\widehat{\gamma})^{(\Omega_{n}^{I}%
,\Omega_{n+1}^{I})})_{\ast}$ we will prove the following theorem.

\begin{theorem}
If either\ (i) $n<p$ or (ii) $n=p\geq1$ and $I=(1,0)$, then the homomorphism
$T(\mathbf{b(}\widehat{\gamma})^{(\Omega_{n}^{I},\Omega_{n+1}^{I})})_{\ast}$
in (1.2) is surjective.
\end{theorem}

These theorem show the importance of the homotopy type of $\Omega^{I}(n,p)$.
In \cite{FoldSurg}, \cite{FoldKyoto}, \cite{FoldRIMS} and \cite{Bo10Topo} we
have determined the homotopy type of $\Omega^{(i,0)}(n,p)$ for $i=\max
\{n-p+1,1\}$, and studied $\mathfrak{NCob}_{n,P}^{(\Omega_{n}^{(i,0)}%
,\Omega_{n+1}^{(i,0)})}$ and $\mathfrak{OCob}_{n,P}^{(\Omega_{n}%
^{(i,0)},\Omega_{n+1}^{(i,0)})}$ in the case $n\geq p$.

When $n=p\geq2$, $I=J=(1,0)$ and $P$ is closed, we will prove in Section 6
that $\pi_{n+\ell}\left(  \widehat{\gamma}_{\boldsymbol{\Omega}_{n+1}^{(1,0)}%
}^{\ell}\right)  $ is, in the oriented case, isomorphic to $\pi_{n+\ell}%
(T(\nu_{P}^{\ell}))$ by using the results in \cite{FoldSurg}. Let $F_{m}$
denote the space of all base point preserving maps of the $m$-sphere $S^{m}$
and let $F=\lim_{m\rightarrow\infty}F_{m}$. By using $S$-dual spaces and
duality maps in the suspension category in \cite{MiSpa} and \cite{SpaDual}, we
can prove that $\pi_{n+\ell}(T(\nu_{P}^{\ell}))$ is isomorphic to the set of
homotopy classes $[P,F]$. Consequently, we have the following theorem.

\begin{theorem}
If $n=p\geq2$ and $P$ is a closed, connected, oriented and $n$-dimensional
manifold, then there exists the isomorphism $\mathfrak{OCob}_{n,P}%
^{(\Omega_{n}^{1,0},\Omega_{n+1}^{1,0})}\rightarrow\lbrack P,F]$.
\end{theorem}

We have constructed the surjection of $\mathfrak{OCob}_{n,P}^{(\Omega
_{n}^{1,0},\Omega_{n+1}^{1,0})}$ onto $[P,F]$ from a different point of view
(\cite[Corollary 2]{FoldKyoto}). This surjection turns out to be a bijection.
Namely, $F$ is the classifying space of this cobordism group. Theorem 1.3
suggests that the stable homotopy groups representing these kind of cobordism
groups yield many invariants related to the singularities of $\Omega^{I}%
$-regular maps.

In \cite{R-S} the group of the cobordism classes of smooth maps of
$n$-dimensional manifolds into $P$ having only a given class of $C^{\infty}$
simple singularities has been represented by the homotopy classes of $P$ to a
certain space in the case $n<p$. The approach leading to their results is
quite different from ours using homotopy principles in this paper. This should
be compared with our Corollary 5.9.

We will actually work in a more general situation than above. We will
generalize the definition of the maps in Theorem 1.1 in Section 3, and will
prove the generalized forms of Theorems 1.1 and 1.2 in Sections 4 and 5
respectively (see Theorem 5.6). Consequently we can apply these theorems to
the groups of cobordism classes of smooth maps having only singularities of
certain $C^{\infty}$ simple types\ by using \cite{hPrinADE} and
\cite{duPContIn}. In Section 6 we will prove Theorem 1.3. We will argue by
applying Theorem 1.3 that stable maps of spheres are detected by higher
singularities through $\mathfrak{OCob}_{n,S^{n}}^{(\Omega_{n}^{1,0}%
,\Omega_{n+1}^{1,0})}$ in low dimensions.

\section{Preliminaries}

Throughout the paper all manifolds are smooth of class $C^{\infty}$. Maps are
basically continuous, but may be smooth (of class $C^{\infty}$) if necessary.
Given a fiber bundle $\pi^{E}:E\rightarrow X$ and a subset $C$ in $X,$ we
denote $\pi^{-1}(C)$ by $E|_{C}.$ Let $\pi^{F}:F\rightarrow Y$ be another
fiber bundle. A map $\tilde{b}:E\rightarrow F$ is called a fiber map over a
map $b:X\rightarrow Y$ if $\pi^{F}\circ\tilde{b}=b\circ\pi^{E}$ holds. The
restriction $\tilde{b}|(E|_{C}):E|_{C}\rightarrow F$ (or $F|_{b(C)}$) is
denoted by $\tilde{b}|_{C}$. In particular, for a point $x\in X,$ $E|_{x}$ and
$\tilde{b}|_{x}$ are simply denoted by $E_{x}$ and $\tilde{b}_{x}%
:E_{x}\rightarrow F_{b(x)}$ respectively. When $E$ and $F$ are vector bundles,
a fiberwise homomorphism, epimorphism and monomorphism $E\rightarrow F$ are
simply called homomorphism, epimorphism and monomorphism respectively. The
trivial bundle $X\times\mathbb{R}^{\ell}$ is denoted by $\varepsilon_{X}%
^{\ell}$.

Let $E\rightarrow X$ (resp. $F\rightarrow Y$) be $n$-dimensional (resp.
$p$-dimensional) vector bundle. Let us recall the $k$-jet bundle $J^{k}%
(E,F)$\ with fiber $J^{k}(n,p)$\ in Introduction, where $k$ may be $\infty$.
The origin of $\mathbb{R}^{m}$ is\ simply denoted by $0$. Let $L^{k}(m)$
denote the group of all $k$-jets of local diffeomorphisms of $(\mathbb{R}%
^{m},0)$. Let $h_{1}:(\mathbb{R}^{p},0)\rightarrow(\mathbb{R}^{p},0)$ and
$h_{2}:(\mathbb{R}^{n},0)\rightarrow(\mathbb{R}^{n},0)$ be local
diffeomorphisms. We define the action of $L^{k}(p)$ $\times$ $L^{k}(n)$ on
$J^{k}(n,p)$ by $(j_{0}^{k}h_{1},j_{0}^{k}h_{2})\cdot j_{0}^{k}f=j_{0}%
^{k}(h_{1}\circ f\circ h_{2}^{-1})$. In particular, $O(p)\times O(n)$ acts on
$J^{k}(n,p)$. Let $\Omega(n,p)$ be an open subset of $J^{k}(n,p)$, which is
invariant with respect to the action of $L^{k}(p)$ $\times$ $L^{k}(n)$. Let
$\Omega(E,F)$ be the open subbundle of $J^{k}(E,F)$ associated to
$\Omega(n,p)$.\ Let $E_{1}\rightarrow X_{1}$ and $F_{1}\rightarrow Y_{1}$ be
other $n$-dimensional vector bundles, and let $\tilde{b}_{1}:E\rightarrow
E_{1}$ and $\tilde{b}_{2}:F\rightarrow F_{1}$ be bundle maps covering
$b_{1}:X\rightarrow X_{1}$ and $b_{2}:Y\rightarrow Y_{1}$ respectively. Then
$\tilde{b}_{1}$ and $\tilde{b}_{2}$\ yield the isomorphisms $S^{i}%
(E_{x})\rightarrow S^{i}(E_{1,b_{1}(x)})$ and $S^{i}(F_{y})\rightarrow
S^{i}(F_{1,b_{2}(y)})$\ for any $x\in X$ and $y\in Y$ for $1\leq i\leq k$ and
hence, we have the bundle map
\begin{equation}
\mathbf{j(}\tilde{b}_{1},\tilde{b}_{2}):J^{k}(E,F)\rightarrow J^{k}%
(E_{1},F_{1})
\end{equation}
covering $b_{1}\times b_{2}$. Then $\mathbf{j(}\tilde{b}_{1},\tilde{b}_{2})$
induces the bundle map $\mathbf{j(}\tilde{b}_{1},\tilde{b}_{2})_{\Omega
}:\Omega(E,F)\rightarrow\Omega(E_{1},F_{1})$.

If we provide $N$ and $P$ with Riemannian metrics, then the Levi-Civita
connections induce the exponential maps $\exp_{N,x}:T_{x}N\rightarrow N$ and
$\exp_{P,y}:T_{y}P\rightarrow P$. In dealing with the exponential maps we
always consider the convex neighborhoods (\cite{KoNo}). We define the smooth
bundle map
\begin{equation}
J^{k}(N,P)\mathbf{\rightarrow}J^{k}(TN,TP)\text{ \ \ \ over }N\times P
\end{equation}
by sending $z=j_{x}^{k}f\in(\pi_{N}^{k}\times\pi_{P}^{k})^{-1}(x,y)$ to the
$k$-jet of $(\exp_{P,y})^{-1}\circ f\circ\exp_{N,x}$ at $\mathbf{0}\in T_{x}%
N$, which is regarded as an element of $J^{k}(T_{x}N,T_{y}P)(=J_{x,y}%
^{k}(TN,TP))$ (see \cite[Proposition 8.1]{KoNo} for the smoothness of
exponential maps). More strictly, (2.2) gives a smooth equivalence of the
fiber bundles under the structure group $L^{k}(p)\times L^{k}(n)$. Namely, it
gives a smooth reduction of the structure group $L^{k}(p)\times L^{k}(n)$ of
$J^{k}(N,P)$ to $O(p)\times O(n)$, which is the structure group of
$J^{k}(TN,TP)$. Let us recall Boardman submanifolds (see \cite{Board},
\cite{LevSDM}). The Boardman submanifold $\Sigma^{I}(N,P)$ of $J^{k}(N,P)$ is
identified with $\Sigma^{I}(TN,TP)$ under (2.2). The same is true for
$\Omega^{I}(N,P)$ and $\Omega^{I}(TN,TP)$.

\section{Maps $\omega_{\mathfrak{N}}^{(\Omega_{n}^{I},\Omega_{n+1}^{J})}$ and
$\omega_{\mathfrak{O}}^{(\Omega_{n}^{I},\Omega_{n+1}^{J})}$}

Let $i_{+1}:J^{k}(n,p)\rightarrow J^{k}(n+1,p+1)$ be the map defined by
$i_{+1}(j_{0}^{k}f)=j_{0}^{k}(f\times id_{\mathbb{R}})$. Let $\Omega(n,p)$ be
given as in Section 2. Let $\Omega^{\star}(n+1,p+1)$ be an open subset of
$J^{k}(n+1,p+1)$ which is invariant with respect to the action of $L^{k}(p+1)$
$\times$ $L^{k}(n+1)$ such that $i_{+1}(\Omega(n,p))\subset\Omega^{\star
}(n+1,p+1)$. Let $\Omega(\gamma_{G_{n}}^{n},TP)$ and $\Omega^{\star}%
(\gamma_{G_{n}}^{n}\oplus\varepsilon_{G_{n}}^{1},TP\oplus\varepsilon_{P}^{1})$
be the open subbundles of $J^{k}(\gamma_{G_{n}}^{n},TP)$ and $J^{k}%
(\gamma_{G_{n}}^{n}\oplus\varepsilon_{G_{n}}^{1},TP\oplus\varepsilon_{P}^{1})$
associated to $\Omega(n,p)$ and $\Omega^{\star}(n+1,p+1)$ respectively. Then
we have the fiber map $\mathbf{i}_{+1}^{(\Omega,\Omega^{\star})}:\Omega
(\gamma_{G_{n}}^{n},TP)\rightarrow\Omega^{\star}(\gamma_{G_{n}}^{n}%
\oplus\varepsilon_{G_{n}}^{1},TP\oplus\varepsilon_{P}^{1})$ over $G_{n}\times
P$ associated to $i_{+1}|\Omega(n,p):\Omega(n,p)\rightarrow\Omega^{\star
}(n+1,p+1)$. Let $\mathbf{c}_{\gamma_{G_{n}}^{n}\oplus\varepsilon_{G_{n}}^{1}%
}:\gamma_{G_{n}}^{n}\oplus\varepsilon_{G_{n}}^{1}\rightarrow\gamma_{G_{n+1}%
}^{n+1}$ be the bundle map, which is induced from, and is covered by
$i^{G}:G_{n}\rightarrow G_{n+1}$. Let us define
\[
\mathbf{j}^{(\Omega,\Omega^{\star})}:\Omega(\gamma_{G_{n}}^{n}%
,TP)\longrightarrow\Omega^{\star}(\gamma_{G_{n+1}}^{n+1},TP\oplus
\varepsilon_{P}^{1})
\]
to be the composite of $\mathbf{i}_{+1}^{(\Omega,\Omega^{\star})}$ and%
\[
\mathbf{j(c}_{\gamma_{G_{n}}^{n}\oplus\varepsilon_{G_{n}}^{1}},id_{TP\oplus
\varepsilon_{P}^{1}})_{\Omega^{\star}}:\Omega^{\star}(\gamma_{G_{n}}^{n}%
\oplus\varepsilon_{G_{n}}^{1},TP\oplus\varepsilon_{P}^{1})\longrightarrow
\Omega^{\star}(\gamma_{G_{n+1}}^{n+1},TP\oplus\varepsilon_{P}^{1}).
\]
We set $\boldsymbol{\Omega}=\Omega(\gamma_{G_{n}}^{n},TP)$,
$\boldsymbol{\Omega}^{\star}=\Omega^{\star}(\gamma_{G_{n+1}}^{n+1}%
,TP\oplus\varepsilon_{P}^{1})$, $\widehat{\gamma}_{\boldsymbol{\Omega}}^{\ell
}=(\pi_{G_{n}}^{k})^{\ast}(\widehat{\gamma}_{G_{n}}^{\ell}%
)|_{\boldsymbol{\Omega}}$ and $\widehat{\gamma}_{\boldsymbol{\Omega}^{\star}%
}^{\ell}=(\pi_{G_{n+1}}^{k})^{\ast}(\widehat{\gamma}_{G_{n+1}}^{\ell
})|_{\boldsymbol{\Omega}^{\star}}$ for simplicity. Since $(i^{G})^{\ast
}(\gamma_{G_{n+1}}^{n+1})=\gamma_{G_{n}}^{n}\oplus\varepsilon_{G_{n}}^{1}$ and
$(i^{G})^{\ast}(\widehat{\gamma}_{G_{n+1}}^{\ell})=\widehat{\gamma}_{G_{n}%
}^{\ell}$, we have the bundle maps%
\begin{align*}
\mathbf{b(}\gamma\oplus\varepsilon^{1})^{(\Omega,\Omega^{\star})}  &
:(\pi_{G_{n}}^{k})^{\ast}(\gamma_{G_{n}}^{n}\oplus\varepsilon_{G_{n}}%
^{1})|_{\boldsymbol{\Omega}}\longrightarrow(\pi_{G_{n+1}}^{k})^{\ast}%
(\gamma_{G_{n+1}}^{n+1})|_{\boldsymbol{\Omega}^{\star}},\\
\mathbf{b(}\widehat{\gamma})^{(\Omega,\Omega^{\star})}  &  :\widehat{\gamma
}_{\boldsymbol{\Omega}}^{\ell}\longrightarrow\widehat{\gamma}%
_{\boldsymbol{\Omega}^{\star}}^{\ell}%
\end{align*}
and the Thom map $T(\mathbf{b(}\widehat{\gamma})^{(\Omega,\Omega^{\star})})$.
Thus we have the homomorphism
\[
T(\mathbf{b(}\widehat{\gamma})^{(\Omega,\Omega^{\star})})_{\ast}:\pi_{n+\ell
}\left(  T(\widehat{\gamma}_{\boldsymbol{\Omega}}^{\ell})\right)
\longrightarrow\pi_{n+\ell}\left(  T(\widehat{\gamma}_{\boldsymbol{\Omega
}^{\star}}^{\ell})\right)  .
\]
We denote the image of $T(\mathbf{b(}\widehat{\gamma})^{(\Omega^{{}}%
,\Omega^{\star})})_{\ast}$ by $\mathfrak{\operatorname{Im}}^{\mathfrak{N}%
}\left(  T(\mathbf{b(}\widehat{\gamma})^{(\Omega,\Omega^{\star})})\right)  $
in the case of $G_{m}=G_{m,\ell}$ or by $\mathfrak{\operatorname{Im}%
}^{\mathfrak{O}}\left(  T(\mathbf{b(}\widehat{\gamma})^{(\Omega,\Omega^{\star
})})\right)  $ in the case of $G_{m}=\widetilde{G}_{m,\ell}$ and $P$ being
oriented respectively.

\begin{definition}
We define $\mathfrak{{NCob}}_{n,P}^{(\Omega,\Omega^{\star})}$\ (resp.
$\mathfrak{OCob}_{n,P}^{(\Omega,\Omega^{\star})}$) to be the set of all (resp.
oriented) $\Omega^{\star}$-cobordism classes of $\Omega$-regular maps of
(resp. oriented) $n$-dimensional manifolds into a connected (resp. oriented)
manifold $P$ by following the definition of $\mathfrak{NCob}_{n,P}%
^{(\Omega_{n}^{I},\Omega_{n+1}^{J})}$\ (resp. $\mathfrak{OCob}_{n,P}%
^{(\Omega_{n}^{I},\Omega_{n+1}^{J})}$) in Introduction and by replacing
$\Omega_{n}^{I}$ and $\Omega_{n+1}^{J}$ by $\Omega$ and $\Omega^{\star}$ respectively.
\end{definition}

Let $\mathfrak{Cob}_{n,P}^{(\Omega,\Omega^{\star})}$ refer to $\mathfrak{NCob}%
_{n,P}^{(\Omega,\Omega^{\star})}$ or $\mathfrak{OCob}_{n,P}^{(\Omega
,\Omega^{\star})}$ and let $\mathrm{\operatorname{Im}}\left(  T(\mathbf{b(}%
\widehat{\gamma})^{(\Omega,\Omega^{\star})})\right)  $ refer to
$\mathrm{\operatorname{Im}}^{\mathfrak{N}}\left(  T(\mathbf{b(}\widehat
{\gamma})^{(\Omega,\Omega^{\star})})\right)  $ or $\mathrm{Im}%
\mathfrak{^{\mathfrak{O}}}\left(  T(\mathbf{b(}\widehat{\gamma})^{(\Omega
,\Omega^{\star})})\right)  $ for simplicity, depending on whether we work in
the nonoriented case or oriented case respectively.

Let $M$ be an $m$-dimensional compact manifold such that $M$ should be
oriented in the oriented case. Take an embedding $e_{M}:M\rightarrow
S^{m+\ell}$ and identify $M$ with $e_{M}(M)$. Let $c_{M}:M\rightarrow G_{m}$
be the classifying map defined by sending a point $x\in M$ to the $m$-plane
$T_{x}M\in G_{m}$. Let $\nu_{M}^{\ell}$\ be the orthogonal normal bundle of
$M$ in $S^{m+\ell}$. Let $\mathbf{c}_{TM}:TM\rightarrow\gamma_{G_{m}}^{m}$
(resp. $\mathbf{c}_{\nu_{M}}:\nu_{M}^{\ell}\rightarrow\widehat{\gamma}_{G_{m}%
}^{\ell}$) be the bundle map covering the classifying map $c_{M}:M\rightarrow
G_{m}$, which is defined by sending a vector $\mathbf{v}$ of $T_{x}M$ (resp.
$\mathbf{w\in}\nu_{M,x}^{\ell}$) to $(T_{x}M,\mathbf{v)}$ (resp.
$(T_{x}M,\mathbf{w)}$). Then we have the canonical trivializations
$t_{M}:TM\oplus\nu_{M}^{\ell}\rightarrow\varepsilon_{M}^{m+\ell}$ and
$t_{G_{m}}:\gamma_{G_{m}}^{m}\oplus\widehat{\gamma}_{G_{m}}^{\ell}%
\rightarrow\varepsilon_{G_{m}}^{m+\ell}$. Then we have $t_{G_{m}}%
\circ(\mathbf{c}_{TM}\oplus\mathbf{c}_{\nu_{M}})\circ t_{M}^{-1}=c_{M}\times
id_{\mathbb{R}^{m+\ell}}$. Furthermore, if there is a map $s_{M}%
:M\rightarrow\Omega(\gamma_{G_{m}}^{m},TP)$ with $\pi_{G_{m}}^{k}\circ s_{M}$
is homotopic (resp. equal) to $c_{M}$, then $\mathbf{c}_{TM}$, $\mathbf{c}%
_{\nu_{M}}$ and the projection $\pi_{G_{m}}^{k}:\Omega(\gamma_{G_{m}}%
^{m},TP)\rightarrow G_{m}$ induce the bundle maps $\mathbf{C}_{TM}%
:TM\rightarrow(\pi_{G_{m}}^{k})^{\ast}\gamma_{G_{m}}^{m}|_{\boldsymbol{\Omega
}}$ and $\mathbf{C}_{\nu_{M}}:\nu_{M}^{\ell}\rightarrow(\pi_{G_{m}}^{k}%
)^{\ast}\widehat{\gamma}_{G_{m}}^{\ell}|_{\boldsymbol{\Omega}}$ covering
$s_{M}$ such that $t_{\mathbf{\Omega}_{m}}\circ(\mathbf{C}_{TM}\oplus
\mathbf{C}_{\nu_{M}})\circ t_{M}^{-1}$ is homotopic (resp. equal) to
$s_{M}\times id_{\mathbb{R}^{m+\ell}}$, where $\mathbf{\Omega}_{m}%
\mathbf{=}\Omega(\gamma_{G_{m}}^{m},TP)$ and $t_{\mathbf{\Omega}_{m}}%
:(\pi_{G_{m}}^{k})^{\ast}(\gamma_{G_{m}}^{m}\oplus\widehat{\gamma}_{G_{m}%
}^{\ell})|_{\mathbf{\Omega}_{m}}\rightarrow\varepsilon_{\mathbf{\Omega}_{m}%
}^{m+\ell}$ is the trivialization induced from $t_{G_{m}}$.

Take an embedding $e_{N}:N\rightarrow S^{n+\ell}$ and apply the above
notation. Then we have the bundle map%
\[
\mathbf{j(c}_{TN},id_{TP})_{\Omega}:\Omega(TN,TP)\rightarrow\boldsymbol{\Omega
}=\Omega(\gamma_{G_{n}}^{n},TP).
\]
Let $f:N\rightarrow P$ be an $\Omega$-regular map and $j^{k}f:N\rightarrow
\Omega(TN,TP)$ be the jet extension of $f$. Then we have the composites%
\[
\mathbf{j}^{(\Omega,\Omega^{\star})}\circ\mathbf{j(c}_{TN},id_{TP})_{\Omega
}\circ j^{k}f:N\rightarrow\boldsymbol{\Omega}^{\star}=\Omega^{\star}%
(\gamma_{G_{n+1}}^{n+1},TP\oplus\varepsilon_{P}^{1}).
\]
and%
\[
\mathbf{b(}\widehat{\gamma})^{(\Omega,\Omega^{\star})}\circ\mathbf{C}_{\nu
_{N}}:\nu_{N}^{\ell}\rightarrow\widehat{\gamma}_{\boldsymbol{\Omega}^{\star}%
}^{\ell}%
\]
covering $\mathbf{j}^{(\Omega,\Omega^{\star})}\circ\mathbf{j(c}_{TN}%
,id_{TP})_{\Omega}\circ j^{k}f$. Let $a_{N}:S^{n+\ell}\rightarrow T(\nu
_{N}^{\ell})$ be the Pontrjagin-Thom construction. We now define the map%
\[
\omega:\mathfrak{Cob}_{n,P}^{(\Omega,\Omega^{\star})}\rightarrow
\mathrm{\operatorname{Im}}\left(  T(\mathbf{b(}\widehat{\gamma})^{(\Omega
,\Omega^{\star})})\right)
\]
by mapping the cobordism class $[f]$ to the homotopy class of $T(\mathbf{b(}%
\widehat{\gamma})^{(\Omega,\Omega^{\star})})\circ T(\mathbf{C}_{\nu_{N}})\circ
a_{N}$. Here, $\omega$ refers to $\omega_{\mathfrak{N}}^{(\Omega,\Omega
^{\star})}$ or $\omega_{\mathfrak{O}}^{(\Omega,\Omega^{\star})}$ depending on
whether we work in the nonoriented case or oriented case. In the rest of the
paper we often deal with the nonoriented case and oriented case at the same
time. We have to prove that $\omega([f])$ does not depend on the choice of a
representative $f$.

\begin{lemma}
Let two $\Omega$-regular maps $f_{i}:N_{i}\rightarrow P$ $(i=0,1)$ are
$\Omega^{\star}$-cobordant. Then we have $\omega_{\mathfrak{N}}^{(\Omega
,\Omega^{\star})}([f_{0}])=\omega_{\mathfrak{N}}^{(\Omega,\Omega^{\star}%
)}([f_{1}])$. If $N_{i}$, $P$ are oriented and $f_{i}$ $(i=0,1)$ are oriented
$\Omega^{\star}$-cobordant, then we have $\omega_{\mathfrak{O}}^{(\Omega
,\Omega^{\star})}([f_{0}])=\omega_{\mathfrak{O}}^{(\Omega,\Omega^{\star}%
)}([f_{1}])$.
\end{lemma}

\begin{proof}
Let $F:(W,\partial W)\rightarrow(P\times\lbrack0,1],P\times0\cup P\times1)$ be
an $\Omega^{\star}$-cobordism of $f_{0}$ and $f_{1}$ as given in Definition
3.1. Take embeddings $e_{N_{i}}:N_{i}\rightarrow S^{n+\ell}$ and
$e_{W}:W\rightarrow S^{n+\ell}\times\lbrack0,1]$ such that

(i) $e_{W}|N_{0}\times\lbrack0,\varepsilon]=e_{N_{0}}\times id_{[0,\varepsilon
]}$,

(ii) $e_{W}|N_{1}\times\lbrack1-\varepsilon,1]=e_{N_{1}}\times
id_{[1-\varepsilon,1]}$.

\noindent We identify $N_{i}$ and $W$ with $e_{N_{i}}(N_{i})$ and $e_{W}(W)$
respectively as above. Let us identify as $TW|_{N_{i}}=TN_{i}\oplus
\varepsilon_{N_{i}}^{1}$. Then we may assume that the trivializations%
\[
t_{N_{i}}:TN_{i}\oplus\ \nu_{N_{i}}^{\ell}\rightarrow\varepsilon_{N_{i}%
}^{n+\ell}\text{ \ \ and \ \ }t_{W}:TW\oplus\ \nu_{W}^{\ell}\rightarrow
\varepsilon_{W}^{n+\ell+1}%
\]
satisfy $t_{W}|_{N_{i}}=(t_{N_{i}}\oplus id_{\varepsilon_{N_{i}}^{1}}%
)\circ(id_{TN_{i}}\oplus\mathbf{k}_{N_{i}}^{\backsim})$, where $\mathbf{k}%
_{N_{i}}^{\backsim}:\varepsilon_{N_{i}}^{1}\oplus\nu_{N_{i}}^{\ell}%
\rightarrow\nu_{N_{i}}^{\ell}\oplus\varepsilon_{N_{i}}^{1}$ is the map
exchanging the components of $\varepsilon_{N_{i}}^{1}$ and $\nu_{N_{i}}^{\ell
}$ and $\varepsilon_{N_{i}}^{n+\ell}\oplus\varepsilon_{N_{i}}^{1}%
=\varepsilon_{N_{i}}^{n+\ell+1}$. Therefore, we have that%
\[%
\begin{array}
[c]{l}%
\mathbf{c}_{\gamma_{G_{n}}^{n}\oplus\varepsilon_{G_{n}}^{1}}\circ
(\mathbf{c}_{TN_{i}}\oplus((i^{G}\circ c_{N_{i}})\times id_{\mathbb{R}%
}))=\mathbf{c}_{TW}|TN_{i}\oplus\varepsilon_{N_{i}}^{1},\\
\mathbf{C}_{\nu_{W}}|_{N_{i}}=\mathbf{b(}\widehat{\gamma})^{(\Omega
,\Omega^{\star})}\circ\mathbf{C}_{\nu_{N_{i}}}.
\end{array}
\]
Let $a_{W}:S^{n+k+1}\times\lbrack0,1]\rightarrow T(\nu_{W}^{\ell})$ be the
Pontrjagin-Thom construction for $e_{W}$. Then the composite $T(\mathbf{C}%
_{\nu_{W}})\circ a_{W}$ gives a homotopy between $T(\mathbf{b(}\widehat
{\gamma})^{(\Omega,\Omega^{\star})})\circ T(\mathbf{C}_{\nu_{N_{0}}})\circ
a_{N_{0}}$ and $T(\mathbf{b(}\widehat{\gamma})^{(\Omega,\Omega^{\star})})\circ
T(\mathbf{C}_{\nu_{N_{1}}})\circ a_{N_{1}}.$ This proves $\omega
([f_{0}])=\omega([f_{1}])$. In the oriented case we only need to provide
manifolds, which appear in the proof, with the orientations.
\end{proof}

\section{$\omega_{\mathfrak{N}}^{(\Omega,\Omega^{\star})}$ and $\omega
_{\mathfrak{O}}^{(\Omega,\Omega^{\star})}$ are isomorphisms}

Let $\mathcal{O}(m,q)$ be an open subset of $J^{k}(m,q)$, which is invariant
with respect to the action $L^{k}(q)\times L^{k}(m)$. Given an $m$-dimensional
manifold $M$ with $\partial M$ and an $q$-dimensional manifold $Q$, we can
define the open subbundle $\mathcal{O}(M,Q)$ in $J^{k}(M,Q)$. Let
$C_{\mathcal{O}}^{\infty}(M,Q)$ denote the space consisting of all
$\mathcal{O}$-regular maps equipped with the $C^{\infty}$-topology. Let
$\Gamma_{\mathcal{O}}(M,Q)$ denote the space consisting of all continuous
sections of the fiber bundle $\pi_{M}^{k}|\mathcal{O}(M,Q):\mathcal{O}%
(M,Q)\rightarrow M$ equipped with the compact-open topology. Then there exists
the continuous map
\[
j_{\mathcal{O}}:C_{\mathcal{O}}^{\infty}(M,Q)\rightarrow\Gamma_{\mathcal{O}%
}(M,Q)
\]
defined by $j_{\mathcal{O}}(f)=j^{k}f$.

\begin{definition}
In this paper we say that $\mathcal{O}(m,q)$ satisfies \textit{the relative
homotopy principle in the existence level }if the following property holds.
Let $M$ and $Q$ be any manifolds as above. Let $C$ be a closed subset of $M$
such that if $\partial M\neq\emptyset$, then $\partial M\subset C$. Let $s$ be
a section of $\Gamma_{\mathcal{O}}(M,Q)$ which has an $\mathcal{O}$-regular
map $g$ defined on a neighborhood of $C$ into $Q$, where $j^{k}g=s$. Then
there exists an $\mathcal{O}$-regular map $f:M\rightarrow Q$ such that
$j^{k}f$ is homotopic to $s$ relative to a neighborhood of $C$ by a homotopy
$s_{\lambda}$ in $\Gamma_{\mathcal{O}}(M,Q)$ with $s_{0}=s$ and $s_{1}=j^{k}f$.
\end{definition}

Let $\Omega(n,p)$ and $\Omega^{\star}(n+1,p+1)$ be an open subset of
$J^{k}(n,p)$ and $J^{k}(n+1,p+1)$ respectively. We say that the pair
$(\Omega(n,p),\Omega^{\star}(n+1,p+1))$ is \textit{admissible to the
h-Principle }if the following properties are satisfied:

(i) $\Omega(n,p)$ (resp. $\Omega^{\star}(n+1,p+1)$) is invariant with respect
to the action of $L^{k}(p)\times L^{k}(n)$ (resp. $L^{k}(p+1)\times
L^{k}(n+1)$).

(ii) $i_{+1}(\Omega(n,p))\subset\Omega^{\star}(n+1,p+1)$.

(iii) $\Omega(n,p)$ and $\Omega^{\star}(n+1,p+1)$ satisfy the relative
homotopy principle in the existence level in Definition 4.1 respectively.

\begin{theorem}
Let $P$ be a connected $p$-dimensional manifold and be oriented in the
oriented case. Assume that the pair $(\Omega(n,p),\Omega^{\star}(n+1,p+1))$ is
\textit{admissible to the h-Principle. Then the maps }%
\begin{align*}
\omega_{\mathfrak{N}}^{(\Omega,\Omega^{\star})}  &  :\mathfrak{NCob}%
_{n,P}^{(\Omega,\Omega^{\star})}\longrightarrow\mathrm{\operatorname{Im}%
}^{\mathfrak{N}}\left(  T(\mathbf{b(}\widehat{\gamma})^{(\Omega,\Omega^{\star
})})\right) \\
\omega_{\mathfrak{O}}^{(\Omega,\Omega^{\star})}  &  :\mathfrak{OCob}%
_{n,P}^{(\Omega,\Omega^{\star})}\longrightarrow\mathrm{\operatorname{Im}%
}^{\mathfrak{O}}\left(  T(\mathbf{b(}\widehat{\gamma})^{(\Omega,\Omega^{\star
})})\right)
\end{align*}
are bijective.
\end{theorem}

\begin{proof}
In the oriented case we only need to provide manifolds which appear in the
proof with the orientations.

We first prove that $\omega$ is injective. For this, take two $\Omega$-regular
maps $f_{i}:N_{i}\rightarrow P$ $(i=0,1)$ such that $\omega([f_{0}%
])=\omega([f_{1}])$. Recall the map $T(\mathbf{b(}\widehat{\gamma}%
)^{(\Omega,\Omega^{\star})})\circ T(\mathbf{c}_{\nu_{N_{i}}})\circ a_{N_{i}}$,
which represents $\omega([f_{i}])$. Let $pr_{X}:X\times\lbrack0,1]\rightarrow
X$ be the canonical projection for a space $X$. There is a homotopy
$H:S^{n+\ell}\times\lbrack0,1]\rightarrow T(\widehat{\gamma}%
_{\boldsymbol{\Omega}^{\star}}^{\ell})$ such that if $\varepsilon$ is
sufficiently small, then

(i) $H|S^{n+\ell}\times\lbrack0,\varepsilon]=T(\mathbf{b(}\widehat{\gamma
})^{(\Omega,\Omega^{\star})})\circ T(\mathbf{c}_{\nu_{N_{0}}})\circ a_{N_{0}%
}\circ(pr_{S^{n+\ell}}|S^{n+\ell}\times\lbrack0,\varepsilon]),$

(ii) $H|S^{n+\ell}\times\lbrack1-\varepsilon,1]=T(\mathbf{b(}\widehat{\gamma
})^{(\Omega,\Omega^{\star})})\circ T(\mathbf{c}_{\nu_{N_{1}}})\circ a_{N_{1}%
}\circ(pr_{S^{n+\ell}}|S^{n+\ell}\times\lbrack1-\varepsilon,1]),$

(iii) $H$ is smooth around $H^{-1}(\boldsymbol{\Omega} ^{\star})$ and is
transverse to $\boldsymbol{\Omega} ^{\star}$.

\noindent We set $W=H^{-1}(\boldsymbol{\Omega} ^{\star})$. Then we have

(iv) $W\cap S^{n+\ell}\times\lbrack0,\varepsilon]=N_{0}\times\lbrack
0,\varepsilon]$ and $H|N_{0}\times\lbrack0,\varepsilon]=\mathbf{j}%
^{(\Omega,\Omega^{\star})}\circ\mathbf{j(c}_{TN_{0}},id_{TP})_{\Omega}\circ
j^{k}f_{0}\circ(pr_{N_{0}}|(N_{0}\times\lbrack0,\varepsilon])$,

(v) $W\cap S^{n+\ell}\times\lbrack1-\varepsilon,1]=N_{1}\times\lbrack
1-\varepsilon,1]$ and $H|N_{1}\times\lbrack1-\varepsilon,1]=\mathbf{j}%
^{(\Omega,\Omega^{\star})}\circ\mathbf{j(c}_{TN_{1}},id_{TP})_{\Omega}\circ
j^{k}f_{1}\circ(pr_{N_{1}}|N_{1}\times\lbrack1-,\varepsilon,1]),$

(vi) $TW|_{N_{0}\times\lbrack0,\varepsilon]}=(TN_{0}\oplus\varepsilon_{N_{0}%
}^{1})\times\lbrack0,\varepsilon]$ and $TW|_{N_{1}\times\lbrack1-\varepsilon
,1]}=(TN_{1}\oplus\varepsilon_{N_{1}}^{1})\times\lbrack1-\varepsilon,1]$.

(vii) $\nu_{W}|_{N_{0}\times\lbrack0,\varepsilon]}=\nu_{N_{0}}\times
\lbrack0,\varepsilon]$ and $\nu_{W}|_{N_{1}\times\lbrack1-\varepsilon,1]}%
=\nu_{N_{1}}\times\lbrack1-\varepsilon,1]$.

By (iii) we have the bundle map $\mathbf{C}_{\nu_{W}}^{\prime}:\nu
_{W}\rightarrow\widehat{\gamma}_{\boldsymbol{\Omega}^{\star}}^{\ell}$ such that

(viii) $\mathbf{C}_{\nu_{W}}^{\prime}|_{N_{0}\times\lbrack0,\varepsilon
]}=\mathbf{C}_{\nu_{N_{0}}}\circ(pr_{\nu_{N_{0}}}|\nu_{N_{0}}\times
\lbrack0,\varepsilon])$ by (i) and (ii),

(ix) $\mathbf{C}_{\nu_{W}}^{\prime}|_{N_{1}\times\lbrack1-\varepsilon
,1]}=\mathbf{C}_{\nu_{N_{1}}}\circ(pr_{\nu_{N_{1}}}|\nu_{N_{1}}\times
\lbrack1-\varepsilon,1])$ by (ii) and (vii).

It follows from \cite[Proposition 3.3]{FoldSurg} that there exists a bundle
map%
\[
\mathbf{C}_{TW}^{\prime}:TW\rightarrow(\pi_{G_{n+1}}^{k})^{\ast}%
(\gamma_{G_{n+1}}^{n+1})|_{\boldsymbol{\Omega}^{\star}}%
\]
covering $H|W:W\rightarrow\boldsymbol{\Omega}^{\star}$ such that
$t_{\boldsymbol{\Omega}^{\star}}\circ(\mathbf{C}_{TW}^{\prime}\oplus
\mathbf{C}_{\nu_{W}}^{\prime})\circ t_{W}^{-1}$ is homotopic to $(H|W)\times
id_{\mathbb{R}^{n+\ell+1}}$. Since $\gamma_{G_{n+1}}^{n+1}$ is the universal
bundle $(\ell\gg n)$, $\mathbf{C}_{TW}^{\prime}$ is homotopic to
$\mathbf{C}_{TW}$. Furthermore, we may assume by (iv), (v) and (vi) that

(x) $\mathbf{C}_{TW}^{\prime}|_{N_{0}\times t}=\mathbf{b(}\gamma
\oplus\varepsilon^{1})^{(\Omega,\Omega^{\star})}\circ(\mathbf{C}%
_{T(N_{0}\times t)}\oplus(c_{_{N_{0}\times t}}\times id_{\mathbb{R}}))$ for
$0\leq t\leq\varepsilon$,

(xi) $\mathbf{C}_{TW}^{\prime}|_{N_{1}\times t}=\mathbf{b(}\gamma
\oplus\varepsilon^{1})^{(\Omega,\Omega^{\star})}\circ(\mathbf{C}%
_{T(N_{1}\times t)}\oplus(c_{_{N_{1}\times t}}\times id_{\mathbb{R}}))$ for
$1-\varepsilon\leq t\leq1$.

\noindent Hence, $c_{W}$ is homotopic to $\pi_{G_{n+1}}^{k}\circ H|W$ relative
to $N_{0}\times\lbrack0,\varepsilon]\cup N_{1}\times\lbrack1-\varepsilon,1]$.
Let $pr_{TP\oplus\varepsilon_{P}^{1}}:T(P\times\lbrack0,1])\rightarrow
TP\oplus\varepsilon_{P}^{1}$ be the canonical bundle map covering
$pr_{P}:P\times\lbrack0,1]\rightarrow P$. Then we have the bundle map%
\[
\mathbf{j(c}_{TW},pr_{TP\oplus\varepsilon_{P}^{1}})_{\mathbf{\Omega}^{\star}%
}:\Omega^{\star}(TW,T(P\times\lbrack0,1]))\longrightarrow\boldsymbol{\Omega
}^{\star}=\Omega^{\star}(\gamma_{G_{n+1}}^{n+1},TP\oplus\varepsilon_{P}^{1})
\]
covering $c_{W}\times pr_{P}$. Therefore, since $[0,1]$ is contractible, there
is the section $s_{W}:W\rightarrow\Omega^{\star}(TW,T(P\times\lbrack0,1])$
such that%
\begin{align*}
\pi_{P\times\lbrack0,1]}^{k}\circ s_{W}|N_{0}\times\lbrack0,\varepsilon]  &
=f_{0}\times id_{[0,\varepsilon]},\\
\pi_{P\times\lbrack0,1]}^{k}\circ s_{W}|N_{1}\times\lbrack1-\varepsilon,1]  &
=f_{1}\times id_{[1-\varepsilon,1],}\\
pr_{P}\circ\pi_{P\times\lbrack0,1]}^{k}\circ s_{W}  &  =\pi_{P}^{k}\circ(H|W),
\end{align*}
and that $\mathbf{j(c}_{TW},pr_{TP\oplus\varepsilon_{P}^{1}})_{\mathbf{\Omega
}^{\star}}\circ s_{W}$ is homotopic to $H|W$ relative to $N_{0}\times
\lbrack0,\varepsilon]\cup N_{1}\times\lbrack1-\varepsilon,1]$.

Since $\Omega^{\star}(TW,T(P\times\lbrack0,1])$ satisfies the relative
homotopy principle in the existence level, there exists an $\Omega^{\star}%
$-regular map $F:W\rightarrow P\times\lbrack0,1]$ such that $F(x,t)=f_{0}%
(x)\times t$ for $0\leq t\leq\varepsilon$, $F(x,t)=f_{1}(x)\times t$ for
$1-\varepsilon\leq t\leq1$ and that $j^{k}F$ is homotopic to $s_{W}$ relative
to $N_{0}\times\lbrack0,\varepsilon/2]\cup N_{1}\times\lbrack1-\varepsilon
/2,1]$. This implies that the $\Omega$-regular maps $f_{0}$ and $f_{1}$ are
$\Omega^{\star}$-cobordant. This proves that $\omega$ is injective.

We next prove that $\omega$ is surjective. Let an element $\widetilde{\alpha}%
$\ of $\mathfrak{\operatorname{Im}}\left(  T(\mathbf{b(}\widehat{\gamma
})^{(\Omega,\Omega^{\star})})\right)  $ be represented by a map $\alpha
:S^{n+\ell}\rightarrow T(\widehat{\gamma}_{\boldsymbol{\Omega}}^{\ell})$ such
that $(T(\mathbf{b(}\widehat{\gamma})^{(\Omega,\Omega^{\star})}))_{\ast
}([\alpha])=\widetilde{\alpha}$. We may suppose that $\alpha$ is smooth around
$\alpha^{-1}(\boldsymbol{\Omega})$ and is transverse to $\boldsymbol{\Omega}$.
we set $N=\alpha^{-1}(\boldsymbol{\Omega})$. If $N=\emptyset$, then $[\alpha]$
must be a null element, while we can deform $\alpha$\ so that $N\neq\emptyset
$\ even in this case.\ Since $\alpha$ is transverse to $\boldsymbol{\Omega}$,
we have the bundle map $\mathbf{C}_{\nu_{N}}^{\prime}:\nu_{N}\rightarrow
\widehat{\gamma}_{\boldsymbol{\Omega}}^{\ell}$ covering $\alpha|N$. It follows
from \cite[Proposition 3.3]{FoldSurg} that there exists a bundle map%
\[
\mathbf{C}_{TN\oplus\varepsilon_{N}^{1}}^{\prime}:TN\oplus\varepsilon_{N}%
^{1}\rightarrow(\pi_{G_{n}}^{k})^{\ast}(\gamma_{G_{n}}^{n}\oplus
\varepsilon_{G_{n}}^{1})|_{\boldsymbol{\Omega}}%
\]
covering $\alpha|N:N\rightarrow\boldsymbol{\Omega}$ such that the composite%
\begin{align*}
&  (t_{\boldsymbol{\Omega}}\oplus id_{\varepsilon_{\boldsymbol{\Omega}}^{1}%
})\circ(id_{(\pi_{G_{n}}^{k})^{\ast}(\gamma_{G_{n}}^{n})}\oplus\mathbf{k}%
_{G_{n}}^{\backsim})\\
&  \text{ \ \ \ \ \ \ \ \ \ \ \ \ \ \ \ \ \ }\circ(\mathbf{C}_{TN\oplus
\varepsilon_{N}^{1}}^{\prime}\oplus\mathbf{C}_{\nu_{N}}^{\prime})\circ
(id_{TN}\oplus\mathbf{k}_{N}^{\backsim})\circ(t_{N}^{-1}\oplus id_{\varepsilon
_{N}^{1}})
\end{align*}
is homotopic to $(\alpha|N)\times id_{\mathbb{R}^{n+\ell+1}}$. Since
$\gamma_{G_{n}}^{n}$ is the universal bundle $(\ell\gg n)$, $\mathbf{C}%
_{TN\oplus\varepsilon_{N}^{1}}^{\prime}$ is homotopic to $\mathbf{C}%
_{TN}\oplus((\alpha|N)\times id_{\mathbb{R}})$, and $t_{\boldsymbol{\Omega}%
}\circ(\mathbf{C}_{TN}\oplus\mathbf{C}_{\nu_{N}}^{\prime})\circ t_{N}^{-1}$ is
homotopic to $(\alpha|N)\times id_{\mathbb{R}^{n+\ell}}$. Hence, $c_{N}$ is
homotopic to $\pi_{G_{n}}^{k}\circ\alpha|N$. By \cite[Proposition
3.3]{FoldSurg} again, $\mathbf{C}_{\nu_{N}}^{\prime}$ and $\mathbf{C}_{\nu
_{N}}$ are homotopic as bundle maps $\nu_{N}\rightarrow\widehat{\gamma
}_{\boldsymbol{\Omega}}^{\ell}$. Since we have the bundle map%
\[
\mathbf{j(c}_{TN},id_{TP})_{\boldsymbol{\Omega}}:\Omega(TN,TP)\longrightarrow
\boldsymbol{\Omega}=\Omega(\gamma_{G_{n}}^{n},TP)
\]
covering $c_{N}$, there is the section $s_{N}:N\rightarrow\Omega(TN,TP)$ such
that $\pi_{P}^{k}\circ s_{N}=\pi_{P}^{k}\circ\alpha|N$ and $\mathbf{j(c}%
_{TN},id_{TP})_{\boldsymbol{\Omega}}\circ s_{N}$ is homotopic to $\alpha|N.$
Since $\Omega(TN,TP)$ satisfies the relative homotopy principle in the
existence level, there exists an $\Omega$-regular map $f:N\rightarrow P$ such
that $j^{k}f$ is homotopic to $s_{N}$. This implies that $\mathbf{j(c}%
_{TN},id_{TP})_{\boldsymbol{\Omega}}\circ j^{k}f$ and $\alpha|N$ are
homotopic. This proves that
\begin{align*}
\omega([f])  &  =[T(\mathbf{b(}\widehat{\gamma})^{(\Omega,\Omega^{\star}%
)})\circ T(\mathbf{C}_{\nu_{N}})\circ a_{N}]\\
&  =[T(\mathbf{b(}\widehat{\gamma})^{(\Omega,\Omega^{\star})})\circ
T(\mathbf{C}_{\nu_{N}}^{\prime})\circ a_{N}]\\
&  =[T(\mathbf{b(}\widehat{\gamma})^{(\Omega,\Omega^{\prime})})\circ\alpha]\\
&  =\widetilde{\alpha}.
\end{align*}
This is what we want.
\end{proof}

Let $\mathcal{O}(m,q)$ be an open subset of $J^{k}(m,q)$, which is invariant
with respect to the action $L^{k}(q)\times L^{k}(m)$. In \cite{duPRS} du
Plessis has called $\mathcal{O}(m,q)$ \textit{extensible} when there exists an
open subset $\mathcal{O}^{\prime}(m+1,q)$ of $J^{k}(m+1,q)$, which is
invariant with respect to the action $L^{k}(q)\times L^{k}(m+1)$, such that
$\widetilde{i}(\mathcal{O}^{\prime}(m+1,q))=\mathcal{O}(m,q)$. Here,
$\widetilde{i}$ is the map which is canonically induced by the inclusion
$i_{\mathbb{R}}:\mathbb{R}^{m}=\mathbb{R}^{m}\times0\rightarrow\mathbb{R}%
^{m+1}$. This \textit{extensibility} yields that not only $j_{\mathcal{O}}$ is
a weak homotopy equivalence, but also $\mathcal{O}(m,q)$ satisfies the
relative homotopy principle in the existence level. However, the last
assertion is not stated explicitly. So we explain an outline of the proof.

\begin{lemma}
Let $\mathcal{O}(m,q)$ be an extensible\textit{ open subset as given above.
Then }$\mathcal{O}(m,q)$ satisfies the relative homotopy principle in the
existence level in Definition 4.1.
\end{lemma}

\begin{proof}
For the $\mathcal{O}$-regular map $g$ and the closed subset $C$ in Definition
4.1, we take a closed neighborhood $V(C)$ of $C$ such that $V(C)$ is an
$n$-dimensional submanifold with boundary and that $g$ is defined on a
neighborhood of $V(C)$, where $j^{k}g=s$. Without loss of generality we may
assume that $N\setminus\mathrm{Int}V(C)$ is nonempty. Take a smooth function
$h_{C}:N\rightarrow\lbrack0,1]$ such that
\begin{equation}
\left\{
\begin{array}
[c]{ll}%
h_{C}(x)=1 & \text{for }x\in C,\\
h_{C}(x)=0 & \text{for }x\in N\setminus\mathrm{Int}V(C),\\
0<h_{C}(x)<1 & \text{for }x\in\mathrm{Int}V(C)\setminus C.
\end{array}
\right.
\end{equation}
By the Sard Theorem (\cite{HirsDT}) there is a regular value $r$ of $h_{C}$
with $0<r<1$. Then $h_{C}^{-1}(r)$ is a submanifold and we set $L_{0}%
=h_{C}^{-1}([r,1])$. We may represent $N$ as the union of an increasing finite
sequence%
\[
L_{0}\subset L_{1}\subset\cdots\subset L_{i}\subset\cdots\subset L_{\sigma}=N
\]
of compact $n$-manifolds with boundary such that $L_{i+1}=L_{i}\cup_{\partial
L_{i}}(\partial L_{i}\times\lbrack0,1])\cup H_{i+1}$, where $H_{i+1}$ is the
$j$-handle $D^{j}\times D^{n-j}$ with $\partial L_{i}\times1\cap
H_{i+1}=\partial D^{j}\times D^{n-j}$. For a sufficiently small $\varepsilon$,
we set $L_{i+1}^{\varepsilon}=L_{i}\cup_{\partial L_{i}}(\partial L_{i}%
\times\lbrack0,\varepsilon])$.

Let $\rho_{C^{\infty}}^{i+1}$ and $\rho_{\Gamma_{\mathcal{O}}}^{i+1}$ in the
following diagram denote the maps, which are canonically induced from the
inclusion $L_{i}^{\varepsilon}\rightarrow L_{i+1}^{\varepsilon}$.%
\[%
\begin{array}
[c]{lll}%
C_{\mathcal{O}}^{\infty}(L_{i+1}^{\varepsilon},P) & \overset{\underrightarrow
{\text{ \ \ }j_{\mathcal{O}}\text{ \ \ }}}{} & \Gamma_{\mathcal{O}}%
(L_{i+1}^{\varepsilon},P)\\
{\scriptstyle{\rho_{C^{\infty}}^{i+1}}}\downarrow &  & \text{ \ \ \ }%
\downarrow{\scriptstyle{\rho_{\Gamma_{\mathcal{O}}}^{i+1}}}\\
C_{\mathcal{O}}^{\infty}(L_{i}^{\varepsilon},P) & \overset{\underrightarrow
{\text{ \ \ }j_{\mathcal{O}}\text{ \ \ }}}{} & \Gamma_{\mathcal{O}}%
(L_{i}^{\varepsilon},P).
\end{array}
\]
It has been proved in \cite{G1} and \cite{duPRS} that $\rho_{C^{\infty}}%
^{i+1}$\ and $\rho_{\Gamma_{\mathcal{O}}}^{i+1}$\ are Serre fibrations. By the
induction on $i\geq0$ we construct an $\mathcal{O}$-regular map $g^{i+1}\in
C^{\infty}(L_{i+1}^{\varepsilon},P)$, a homotopy $u_{\lambda}^{i+1}\in
\Gamma_{\mathcal{O}}(L_{i+1}^{\varepsilon},P)$\ relative to $L_{0}$\ and a
homotopy $s_{\lambda}^{i+1}\in\Gamma_{\mathcal{O}}(L_{i+1}^{\varepsilon}%
,P)$\ relative to $L_{i}$\ such that $u_{0}^{i+1}=s|L_{i+1}^{\varepsilon}$,
$u_{1}^{i+1}|L_{i}=j^{k}g^{i}$, $g^{i+1}|L_{i}=g^{i}$, $u_{1}^{i+1}%
=s_{0}^{i+1}$, $s_{0}^{i+1}|L_{i}=j^{k}g^{i}$ and $s_{1}^{i+1}=j^{k}g^{i+1}$.
Let $v_{\lambda}^{i+1}\in\Gamma_{\mathcal{O}}(L_{i+1}^{\varepsilon},P)$ be the
homotopy defined by $v_{\lambda}^{i+1}=u_{2\lambda}^{i+1}$ ($0\leq\lambda
\leq1/2$) and $v_{\lambda}^{i+1}=u_{2\lambda-1}^{i+1}$ ($1/2\leq\lambda\leq1$).

In fact, we start with $g^{0}=g|L_{0}^{\varepsilon}$ and $u_{\lambda}%
^{0}=s_{\lambda}^{0}=s|L_{0}^{\varepsilon}$\ for any $\lambda,$\ and\ next
assume that $g^{i}$, $u_{\lambda}^{i}$ and $s_{\lambda}^{i}$ are already
constructed. By applying the homotopy extension property to $s|L_{i+1}%
^{\varepsilon}$ and $v_{\lambda}^{i}$, we have the homotopy $u_{\lambda}%
^{i+1}\in\Gamma_{\mathcal{O}}(L_{i+1}^{\varepsilon},P)$\ relative to $L_{0}%
$\ such that $u_{0}^{i+1}=s|L_{i+1}^{\varepsilon}$, $u_{\lambda}^{i+1}%
|L_{i}^{\varepsilon}=v_{\lambda}^{i}|L_{i}^{\varepsilon}$. Let $\mathcal{F}%
_{C^{\infty}}^{i+1}$\ and $\mathcal{F}_{\Gamma_{\mathcal{O}}}^{i+1}$ denote
the fibers of $\rho_{C^{\infty}}^{i+1}$\ and $\rho_{\Gamma_{\mathcal{O}}%
}^{i+1}$\ over $g^{i}$ and $j^{k}g^{i}$ respectively. Since $j_{\mathcal{O}}$
induces the weakly homotopy equivalence $\mathcal{F}_{C^{\infty}}%
^{i+1}\rightarrow\mathcal{F}_{\Gamma_{\mathcal{O}}}^{i+1}$, we have an
$\mathcal{O}$-regular map $g^{i+1}$ and a homotopy $s_{\lambda}^{i+1}%
$\ relative to $L_{i}$\ such that $g^{i+1}|L_{i}=g^{i}$, $s_{0}^{i+1}%
=u_{1}^{i+1}$, $s_{1}^{i+1}=j^{\infty}g^{i+1}$. This is what we want.

Define $s_{\lambda}\in\Gamma_{\mathcal{O}}(L_{\sigma},P)$ to be the homotopy
$v_{\lambda}^{\sigma}$. Then we have $s_{0}=s$ and the required $\mathcal{O}%
$-regular map $f=g^{\sigma}$ with $s_{1}=j^{k}g^{\sigma}$.
\end{proof}

\begin{proof}
[Proof of Theorem 1.1]Let $I=(i_{1},i_{2},\cdots,i_{k})$. It has been proved
in \cite{duPMCS} that if $i_{k}>n-p-d^{I}$, where $d^{I}$ is the sum of
$\alpha_{1},\cdots\alpha_{k-1}$ with $\alpha_{\ell}$ being $1$ or $0$
depending on $i_{\ell}-i_{\ell+1}>1$ or otherwise, then $\Omega^{I}(n,p)$ is extensible.

Let $n<p$.\ By Lemma 4.3, $\Omega^{I}(n,p)$ and $\Omega^{J}(n+1,p+1)$\ satisfy
the relative homotopy principle in the existence level (see also
\cite{HirsImm} and \cite{Feit}).

In \cite[Theorem 0.1]{hPrinBo} it has been proved that if $n\geq p\geq2$ and
$I\geq(n-p+1,0)$, then $\Omega^{I}(n,p)$ and $\Omega^{J}(n+1,p+1)$\ satisfy
the relative homotopy principle in the existence level. Here we remark the
following. In \cite[Theorem 0.1]{hPrinBo} $V$ is assumed to be $\partial
V=\emptyset$. However, this does not matter, because we only need to consider
the manifold $V-\partial V$.

Therefore, the pair $(\Omega^{I}(n,p),\Omega^{J}(n+1,p+1))$ is admissible to
the h-Principle if (i) $n<p$ or (ii) $n\geq p$ and $I\geq(n-p+1,0)$. This
proves the theorem.
\end{proof}

In \cite[Section 0, Theorem]{duPContIn} there have been given extensible open
subsets $\Omega(n,p)$, which are associated to smooth maps having only
singularities of certain $C^{\infty}$ simple type. In \cite{hPrinADE} we have
constructed the submanifolds $\Sigma D_{i}(n,p)$ ($i\geq4$) and $\Sigma
E_{i}(n,p)$ ($i=6$, $7$ or $8$) in $J^{k}(n,p)$, which play the similar role
for the singularities of types $D_{i}$ and $E_{i}$ respectively as
$\Sigma^{(n-p+1,1,\cdots,1,0)}(n,p)$ does for the singularities of type
$A_{i}$ (see \cite{Arnold}). Consider a subset $\Omega(n,p)$ of $J^{k}(n,p)$,
which consists of all regular jets and a number of prescribed submanifolds
$\Sigma^{(n-p+1,1,\cdots,1,0)}(n,p)$, $\Sigma D_{i}(n,p)$ and $\Sigma
E_{j}(n,p)$. We can find when $\Omega(n,p)$ becomes open by using the
adjacency relations of these singularities given in \cite[Corollary
8.7]{Arnold}. We can apply Theorem 4.2 to these open subsets $\Omega(n,p)$
(see also Corollary 5.9).

Let us provide $\mathfrak{Cob}_{n,P}^{(\Omega,\Omega^{\star})}$ with the
structure of a module so that $\omega$\ is an isomorphism. This is a standard
argument and the details are left to the reader. Given two $\Omega$-regular
maps $f_{i}:N_{i}\rightarrow P$ ($i=0,1$), we define the sum $[f_{0}]+[f_{1}]$
to be the cobordism class of the $\Omega$-regular maps $f:N_{0}\cup
N_{1}\rightarrow P$ defined by $f|N_{i}=f_{i}$. The null element is defined to
be represented by an $\Omega$-regular maps $f:N\rightarrow P$, which has an
$\Omega^{\star}$-cobordism $F:(W,\partial W)\rightarrow(P\times\lbrack
0,1],P\times0)$ with $\partial W=N$ such that $F|N=f$ under the identification
$P\times0=P$.

\section{Proof of Theorem 1.2}

In this section we prove Theorem 1.2. Let $C_{m}$ denote the set consisting of
all smooth map germs $(\mathbb{R}^{m},0)\rightarrow\mathbb{R}$ and let
$\mathfrak{m}_{m}$ denote the ideal in $C_{m}$ which consists of all smooth
map germs vanishing at the origin. For a $k$-jet $z=j^{k}f\in J^{k}(m,q)$\ we
define the $\mathbb{R}$-algebra $Q(z)=C_{m}/(f^{\ast}(\mathfrak{m}%
_{q})+\mathfrak{m}_{m}^{k+1})$.\ If two $\mathbb{R}$-algebras $A$%
\ and\ $B$\ are isomorphic, then we write $A\approx B$. Let $\Omega(n,p)$\ and
$\Omega(n+1,p+1)$\ be subsets of\ $J^{k}(n,p)$\ and $J^{k}(n+1,p+1)$%
\ respectively\ satisfying the conditions:

(\textbf{C}1) $\Omega(n,p)$\ and $\Omega(n+1,p+1)$\ are open and
$i_{+1}(\Omega(n,p))\subset\Omega(n+1,p+1)$.

(\textbf{C}2) If a $k$-jet $z\in J^{k}(n,p)$ (resp. $z\in J^{k}(n+1,p+1)$%
)\ has a $k$-jet $w\in\Omega(n,p)$ (resp. $w\in\Omega(n+1,p+1)$) such that
$Q(z)\approx Q(w)$, then $z\in\Omega(n,p)$ (resp. $z\in\Omega(n+1,p+1)$).

(\textbf{C}3) If a $k$-jet $z\in\Omega(n+1,p+1)$\ does not lie in
$\Sigma^{n+1}(n+1,p+1)$, then there exists\ a $k$-jet $w\in\Omega(n,p)$ such
that $Q(z)\approx Q(w)$.

By (\textbf{C}2), $\Omega(n,p)$\ and $\Omega(n+1,p+1)$\ are invariant with
respect to the actions of $L^{k}(p)\times L^{k}(n)$ and $L^{k}(p+1)\times
L^{k}(n+1)$ respectively.

\begin{lemma}
The open subsets $\Omega^{I}(n,p)$ and $\Omega^{I}(n+1,p+1)$ for the symbol
$I=(i_{1},\cdots,i_{k})$ satisfy the conditions (\textbf{C}1), (\textbf{C}2)
and (\textbf{C}3).

\begin{proof}
The condition (\textbf{C}1) and (\textbf{C}2) are satisfied by \cite[2,
Corollary(Morin)]{MathTB}.

We consider the usual coordinates $x=(x_{1},x_{2},\cdots,x_{n+1})$ of
$\mathbb{R}^{n+1}$ and $y=(y_{1},y_{2},\cdots,y_{p+1})$ of $\mathbb{R}^{p+1}$.
Since the given $k$-jet $z\in\Omega(n+1,p+1)$\ does not lie in $\Sigma
^{n+1}(n+1,p+1)$, $z$ is represented as $z=j_{0}^{k}f$ with%
\[
(y_{1}^{\prime}\circ f(x^{\prime}),\cdots,y_{p+1}^{\prime}\circ f(x^{\prime
}))=(g^{1}(x^{\prime}),\cdots,g^{p-n+i_{1}}(x^{\prime}),x_{i_{1}+1}^{\prime
},\cdots,x_{n+1}^{\prime}),
\]
where $g^{j}\in\mathfrak{m}_{n+1}^{2}$ under suitable coordinates $x^{\prime
}=(x_{1}^{\prime},x_{2}^{\prime},\cdots,x_{n+1}^{\prime})$ of $\mathbb{R}%
^{n+1}$ and $y^{\prime}=(y_{1}^{\prime},y_{2}^{\prime},\cdots,y_{p+1}^{\prime
})$ of $\mathbb{R}^{p+1}$. Let $V$ and $W$ be the subspaces of $\mathbb{R}%
^{n+1}$ and $\mathbb{R}^{p+1}$\ defined by the equations $x_{n+1}^{\prime}%
=0$\ and $y_{p+1}^{\prime}=0$, and\ let $i_{V}:V\rightarrow\mathbb{R}^{n+1}%
$\ and $\pi_{W}:\mathbb{R}^{p+1}\rightarrow W$\ be the inclusion and the
projection defined by $\pi_{W}(y_{1}^{\prime},y_{2}^{\prime},\cdots
,y_{p+1}^{\prime})=(y_{1}^{\prime},y_{2}^{\prime},\cdots,y_{p}^{\prime}%
)$\ respectively. Let $\overset{\bullet}{x}=(x_{1},x_{2},\cdots,x_{n})$ and
$\overset{\bullet}{x}^{\prime}=(x_{1}^{\prime},x_{2}^{\prime},\cdots
,x_{n}^{\prime})$. Define the map $\overline{g}:V\rightarrow W$\ by
$\overline{g}=\pi_{W}\circ f\circ i_{V}$ and the functions $\overline{g^{j}%
}\,$on $V$ by $\overline{g^{j}}(\overset{\bullet}{x}^{\prime})=y_{j}^{\prime
}\circ\overline{g}(\overset{\bullet}{x}^{\prime})$\ ($1\leq j\leq p$).
Replacing the coordinates $x_{j}^{\prime}$ and $y_{j}^{\prime}$\ by\ $x_{j}$
and $y_{j}$\ we define the germ $\overline{f}:\mathbb{R}^{n}\rightarrow
\mathbb{R}^{p}$\ by
\[
\overline{f}(\overset{\bullet}{x})=\left\{
\begin{array}
[c]{ll}%
(\overline{g^{1}}(\overset{\bullet}{x}),\overline{g^{2}}(\overset{\bullet}%
{x}),\cdots,\overline{g^{p-n+i_{1}}}(\overset{\bullet}{x}),x_{i_{1}+1}%
,\cdots,x_{n}) & \text{for }i_{1}<n,\\
(\overline{g^{1}}(\overset{\bullet}{x}),\overline{g^{2}}(\overset{\bullet}%
{x}),\cdots,\overline{g^{p}}(\overset{\bullet}{x})) & \text{for }i_{1}=n.
\end{array}
\right.
\]

If we write $g(x^{\prime})$ for $f(x^{\prime})$ to avoid the confusion, then%
\begin{align*}
Q(z)  &  =C_{n+1}/(f^{\ast}(\mathfrak{m}_{p+1})+\mathfrak{m}_{n+1}^{k+1})\\
&  =C_{n+1}/(g^{\ast}(\mathfrak{m}_{p+1})+\mathfrak{m}_{n+1}^{k+1})\\
&  \approx C_{n}/(\overline{g}^{\ast}(\mathfrak{m}_{p})+\mathfrak{m}_{n}%
^{k+1})\\
&  \approx C_{n}/(\overline{f}^{\ast}(\mathfrak{m}_{p})+\mathfrak{m}_{n}%
^{k+1})\\
&  =Q(j_{0}^{k}\overline{f}).
\end{align*}
By \cite[2, Corollary(Morin)]{MathTB}, the Thom-Boardman symbol of $j_{0}%
^{k}\overline{f}$ is equal to $I$. This proves the assertion.
\end{proof}
\end{lemma}

A point of$\ \mathbb{R}^{m}$ is expressed as $x_{1}\mathbf{e}_{1}%
+x_{2}\mathbf{e}_{2}+\cdots+x_{m}\mathbf{e}_{m}$\ or $(x_{1},x_{2}%
,\cdots,x_{m})$, where $\mathbf{e}_{1}$, $\mathbf{e}_{2}$, $\cdots$,
$\mathbf{e}_{m}$\ are the canonical orthonormal basis. Let $pr_{p+1}%
:\mathbb{R}^{p+1}\rightarrow\mathbb{R}$ be the projection mapping $^{t}%
(y_{1},y_{2},\cdots,y_{p+1})$ to $y_{p+1}$. Let $j(pr_{p+1}):J^{1}%
(n+1,p+1)\rightarrow J^{1}(n+1,1)$ denote the map defined by mapping a $1$-jet
$j_{\mathbf{0}}^{1}f$ to the $1$-jet $j_{\mathbf{0}}^{1}(pr_{p+1}\circ f)$.

Let $K$ be a finite simplicial complex of dimension $i\ $and$\ L$ be its
subcomplex of dimension less than $i$ such that $K\backslash L$ is a manifold.

\begin{lemma}
Let $\Omega(n,p)$\ and $\Omega(n+1,p+1)$\ satisfy (\textbf{C}1), (\textbf{C}2)
and (\textbf{C}3) as above. Let $i\leq n<p$ and $(K,L)$ be given as above. Let
$\psi:(K,L)\rightarrow(\Omega(n+1,p+1),i_{+1}(\Omega(n,p)))$ be a map such
that $\psi|(K\backslash L)$ is smooth. Then there exists a homotopy
$\psi_{\lambda}:(K,L)\rightarrow(\Omega(n+1,p+1),i_{+1}(\Omega(n,p)))$ such that

(i) $\psi_{0}=\psi,$

(ii) $\psi_{\lambda}|L=\psi|L,$

(iii) $j(pr_{p+1})\circ\pi_{1}^{k}\circ\psi_{1}(u)(\mathbf{e}_{i})=\left\{
\begin{array}
[c]{ll}%
\mathbf{0} & \text{for }i<n+1,\\
\mathbf{e}_{p+1} & \text{for }i=n+1.
\end{array}
\right.  $
\end{lemma}

\begin{proof}
Let us define $\kappa:K\rightarrow\mathbb{R}^{p+1}$ by $\kappa(u)=(\pi_{1}%
^{k}\circ\psi(u))(\mathbf{e}_{n+1})$. Since $\psi(L)\subset i_{+1}%
(\Omega(n,p))$, we have that, for any $u\in L$, $\kappa(u)=\mathbf{e}_{p+1}$.
Since $\dim K\leq n$, $K\backslash L$ is a manifold and since $\psi
|(K\backslash L)$ is smooth, there exists a deformation $\kappa_{\lambda}$ of
$\kappa$ relative to a neighborhood of $L$ with $\kappa_{0}=\kappa$ such that
the set of $\lambda$'s for which $\kappa_{\lambda}$ does not take the value
$\mathbf{0}\in\mathbb{R}^{p+1}$ is dense in $[0,1]$. Consider the fiber bundle
$q_{\mathbb{R}^{p+1}}:J^{k}(n+1,p+1)\rightarrow\mathbb{R}^{p+1}$ defined by
$q_{\mathbb{R}^{p+1}}(j_{\mathbf{0}}^{k}f)=j_{\mathbf{0}}^{1}f(\mathbf{e}%
_{n+1})$. By applying the covering homotopy property to $\psi$ and
$\kappa_{\lambda}$, there exists a homotopy $\varphi_{\lambda}$ relative to
$L$ such that $\varphi_{0}=\psi$ and $q_{\mathbb{R}^{p+1}}\circ\varphi
_{\lambda}=\kappa_{\lambda}$. Since $\Omega(n+1,p+1)$ is an open subset and
$K$ is compact, there exists a $r\in\lbrack0,1]$ such that if $\lambda\leq r$,
then $\varphi_{\lambda}(K)\subset\Omega(n+1,p+1)$ and $\kappa_{r}$ does not
take the value $\mathbf{0\in}\mathbb{R}^{p+1}$.

In the proof an element of $GL(m)$ is regarded as a linear isomorphism of
$\mathbb{R}^{m}$.\ Let $h_{\lambda}^{1}:(K,L)\rightarrow(GL(p+1),E_{p+1})$ be
the homotopy defined by $h_{\lambda}^{1}(u)=((1-\lambda)+\lambda/\Vert
\kappa_{r}(u)\Vert)E_{p+1}$. Then we have $h_{1}^{1}(u)(\kappa_{r}(u))\in
S^{p}$, and we may assume without loss of generality that $h_{1}^{1}%
(u)(\kappa_{r}(u))\neq-\mathbf{e}_{p+1}$ for any $u\in K$. By considering the
rotation which is the identity on all points orthogonal to both $\kappa
_{r}(u)$ and $\mathbf{e}_{p+1}$ and rotate the great circle through
$\kappa_{r}(u)$ and $\mathbf{e}_{p+1}$ so as to carry $\kappa_{r}(u)$ to
$\mathbf{e}_{p+1}$ (when $\kappa_{r}(u)=\mathbf{e}_{p+1}$, we consider
$E_{p+1}$), we have the homotopy $h_{\lambda}^{2}:(K,L)\rightarrow
(SO(p+1),E_{p+1})$ relative to $L$ such that $h_{0}^{2}(u)=E_{p+1}$ and
$h_{1}^{2}(u)\circ h_{1}^{1}(u)(\kappa_{r}(u))=\mathbf{e}_{p+1}$ for any $u\in
K$. Let $h_{\lambda}:(K,L)\rightarrow(GL(p+1),E_{p+1})$ be the homotopy
defined by $h_{\lambda}=h_{2\lambda}^{1}$ for $0\leq\lambda\leq1/2$ and
$h_{\lambda}=h_{2\lambda-1}^{2}\circ h_{1}^{1}$ for $1/2\leq\lambda\leq1$.
Define $\mathcal{K}^{J}:K\rightarrow J^{1}(n+1,1)$ by
\[
\mathcal{K}^{J}(u)=j(pr_{p+1})\circ\pi_{1}^{k}\circ(h_{1}(u)\circ\varphi
_{r}(u)).
\]
Since $\mathcal{K}^{J}(u)$ is of rank $1$ for any $u\in K$, we have the unique
vector $\mathbf{v(}u)$ of length $1$ such that $\mathbf{v}(u)$ is
perpendicular to $\mathrm{Ker}(\mathcal{K}^{J}(u))$ and that $\mathcal{K}%
^{J}(u)(\mathbf{v}(u))$ is directed to the same orientation of $\mathbf{e}%
_{p+1}$. Since $h_{1}(u)(\kappa_{r}(u))=\mathbf{e}_{p+1}$, $\mathbf{v}(u)$
cannot be equal to $-\mathbf{e}_{n+1}$. We set $v(u)=\Vert\mathcal{K}%
^{J}(u)(\mathbf{v}(u))\Vert$.

Note that $\mathbf{v}(u)\neq-\mathbf{e}_{n+1}$ for any $u\in K$. By
considering the rotation which is the identity on all points orthogonal to
both $\mathbf{v}(u)$ and $\mathbf{e}_{n+1}$ and rotate the great circle
through $\mathbf{v}(u)$ and $\mathbf{e}_{n+1}$ so as to carry $\mathbf{e}%
_{n+1}$ to $\mathbf{v}(u)$, we again have the homotopy $H_{\lambda}%
^{1}:(K,L)\rightarrow(SO(n+1),E_{n+1})$ relative to $L$ such that $H_{0}%
^{1}(u)=E_{n+1}$ and $H_{1}^{1}(u)(\mathbf{e}_{n+1})=\mathbf{v}(u)$ for any
$u\in K$. Let$\ H_{\lambda}^{2}:(K,L)\rightarrow(GL(n+1),E_{n+1})$ be the
homotopy relative to $L$ defined by $H_{\lambda}^{2}(u)=((1-\lambda
)+\lambda/v(u))E_{n+1}$. Let $H_{\lambda}:(K,L)\rightarrow(GL(n+1),E_{n+1})$
be the homotopy defined by $H_{\lambda}(u)=H_{2\lambda}^{1}(u)$ for
$0\leq\lambda\leq1/2$ and $H_{\lambda}(u)=H_{2\lambda-1}^{2}(u)\circ H_{1}%
^{1}(u)$ for $1/2\leq\lambda\leq1$. Then we have that, for any $u\in K$,
\begin{align*}
\mathcal{K}^{J}(u)\circ H_{1}(u)(\mathbf{e}_{n+1})  &  =\mathcal{K}%
^{J}(u)\circ H_{1}^{2}(u)\circ H_{1}^{1}(u)(\mathbf{e}_{n+1})\\
&  =\mathcal{K}^{J}(u)\circ H_{1}^{2}(u)(\mathbf{v}(u))\\
&  =\mathcal{K}^{J}(u)(\mathbf{v}(u))/v(u)\\
&  =\mathbf{e}_{p+1}.
\end{align*}
Since $H_{1}^{1}(u)\in SO(n+1)$ and $\mathbf{e}_{i}$ is orthogonal to
$\mathbf{e}_{n+1}$\ ($i<n+1$), $H_{1}^{1}(u)(\mathbf{e}_{i})$ is orthogonal to
$H_{1}^{1}(u)(\mathbf{e}_{n+1})=\mathbf{v}(u)$. Namely, $H_{1}(u)(\mathbf{e}%
_{i})$ lies in $\mathrm{Ker}(\mathcal{K}^{J}(u))$. Hence, we have%
\[
\mathcal{K}^{J}(u)\circ H_{1}(u)(\mathbf{e}_{i})=\mathbf{0}\text{ \ \ \ for
}i<n+1.
\]
Define the homotopy $\psi_{\lambda}:(K,L)\rightarrow(\Omega(n+1,p+1),i_{+1}%
(\Omega(n,p))$ relative to $L$ by
\[
\psi_{\lambda}(u)=\left\{
\begin{array}
[c]{ll}%
\varphi_{3\lambda r}(u) & \text{for }0\leq\lambda\leq1/3,\\
h_{3\lambda-1}(u)\circ\varphi_{r}(u) & \text{for }1/3\leq\lambda\leq2/3,\\
h_{1}(u)\circ\varphi_{r}(u)\circ H_{3\lambda-2}(u) & \text{for\ }%
2/3\leq\lambda\leq1.
\end{array}
\right.
\]
By the definition we have
\[
j(pr_{p+1})\circ\pi_{1}^{k}\circ\psi_{1}(u)(\mathbf{e}_{i})=\left\{
\begin{array}
[c]{ll}%
\mathbf{0} & \text{for }i<n+1,\\
\mathbf{e}_{p+1} & \text{for }i=n+1.
\end{array}
\right.
\]

This is what we want.
\end{proof}

\begin{proposition}
Under the same assumption of Lemma 5.2, we have a homotopy $\Psi_{\lambda
}:(K,L)\rightarrow(\Omega(n+1,p+1),i_{+1}(\Omega(n,p)))$ such that

(i) $\Psi_{0}=\psi,$

(ii) $\Psi_{\lambda}|L=\psi|L,$

(iii) $\Psi_{1}(K)\subset i_{+1}(\Omega(n,p)).$
\end{proposition}

\begin{proof}
Let $\psi_{\lambda}$ be the homotopy given in Lemma 5.2. Let us express
$\psi_{\lambda}(u)=(f_{\lambda}^{1}(u),f_{\lambda}^{2}(u),\cdots,f_{\lambda
}^{p+1}(u))$\ using the coordinates of $\mathbb{R}^{p+1}$, where $f_{\lambda
}^{i}(u)$\ is regarded as a polynomial of degree $k$ with constant $0$. We
note that%
\[
f_{\lambda}^{p+1}(u)(x_{1},\cdots,x_{n+1})=x_{n+1}+\text{higher term}.
\]
Let $C_{\mathbf{0}}^{\infty}(\mathbb{R}^{n+1},\mathbb{R}^{n+1})$\ denote the
set of all germs of local diffeomorphisms of $(\mathbb{R}^{n+1},\mathbf{0)}$.
Let us define a homotopy of maps $\Phi_{\lambda}:(K,L)\rightarrow
C_{\mathbf{0}}^{\infty}(\mathbb{R}^{n+1},\mathbb{R}^{n+1})$ by%
\[
\Phi_{\lambda}(u)(x_{1},\cdots,x_{n+1})=(x_{1},\cdots,x_{n+1}+\lambda
(f_{\lambda}^{p+1}(u)-x_{n+1})).
\]
It is obvious that $\Phi_{\lambda}(u)$\ is a germ of a diffeomorphism of
$(\mathbb{R}^{n+1},\mathbf{0)}$. Then we have the inverse $\Phi_{\lambda
}(u)^{-1}$ such that%
\begin{equation}
pr_{p+1}\circ\psi_{1}(u)\circ\Phi_{1}(u)^{-1}(x_{1},\cdots,x_{n+1}%
)=x_{n+1}\text{.}%
\end{equation}
We now define $\phi_{\lambda}:(K,L)\rightarrow(\Omega(n+1,p+1),i_{+1}%
(\Omega(n,p)))$ by%
\[
\phi_{\lambda}(u)=\psi_{1}\circ j_{\mathbf{0}}^{k}(\Phi_{\lambda}(u)^{-1}).
\]

Next we exclude the terms containing $x_{n+1}$ from $y_{j}\circ\phi_{\lambda}%
$\ ($1\leq j\leq p$). Let $\eta_{\lambda}:(K,L)\rightarrow(J^{k}%
(n+1,p+1),J^{k}(n,p)))$\ be the homotopy defined by
\begin{align}
\eta_{\lambda}(u)(x)  &  =(1-\lambda)\phi_{1}(u)((x_{1},\cdots,x_{n+1}%
)\nonumber\\
&  +\lambda(\phi_{1}(u)(x_{1},\cdots,x_{n},0)+(0,\cdots,0,x_{n+1})).
\end{align}
It is obvious that $\eta_{1}(K)\subset i_{+1}(\Omega(n,p))$ and that
$\eta_{\lambda}(L)=\psi|L.$ It remains to prove that $\eta_{\lambda}$ is a
homotopy into $\Omega(n+1,p+1)$. It follows from (5.1) and (5.2) that%
\[
pr_{p+1}\circ\eta_{\lambda}(u)(x)=(1-\lambda)(x_{n+1})+\lambda x_{n+1}%
=x_{n+1}.
\]
Let us express $\eta_{\lambda}(u)=(g_{\lambda}^{1}(u),g_{\lambda}%
^{2}(u),\cdots,g_{\lambda}^{p+1}(u))$, where $g_{i}(u)$\ is regarded as a
polynomial of degree $k$ with constant $0$.\ Consider the ideal $\mathfrak{I}%
_{\lambda}(u)$\ generated by $g_{\lambda}^{1}(u),g_{\lambda}^{2}%
(u),\cdots,g_{\lambda}^{p+1}(u)$\ in $\mathfrak{m}_{n+1}\mathfrak{/m}%
_{n+1}^{k+1}$. Then $\mathfrak{I}_{\lambda}(u)$\ is constantly equal to
$\mathfrak{I}_{0}(u)$, and hence $Q(\eta_{\lambda}(u))\approx Q(\psi_{1}(u))$.
Since $\psi_{1}(u)\in\Omega(n+1,p+1)$, we have $\eta_{\lambda}(u)\in
\Omega(n+1,p+1)$ by (\textbf{C}2). Then the required homotopy $\Psi_{\lambda}$
is defined by $\Psi_{\lambda}=\psi_{3\lambda}\,$($0\leq\lambda\leq1/3$),
$\Psi_{\lambda}=\phi_{3\lambda-1}\ $($1/3\leq\lambda\leq2/3$) and
$\Psi_{\lambda}=\eta_{3\lambda-2}$ ($2/3\leq\lambda\leq1$).
\end{proof}

\begin{proposition}
Let $n<p$. Let $\Omega(n,p)$\ and $\Omega(n+1,p+1)$\ denote the open
subspaces, which satisfy (\textbf{C}1), (\textbf{C}2) and (\textbf{C}3). Then
$(i_{+1})_{\ast}:\pi_{i}(\Omega(n,p))\rightarrow\pi_{i}(\Omega(n+1,p+1))$ is
an isomorphism for $0\leq i<n$ and an epimorphism for $i=n$.
\end{proposition}

\begin{proof}
Let $\iota_{n}:\mathbb{R}^{n}\rightarrow\mathbb{R}^{p}$ and $\iota
_{n+1}:\mathbb{R}^{n+1}\rightarrow\mathbb{R}^{p+1}$ denote the map defined by%
\[%
\begin{array}
[c]{l}%
\iota_{n}(x_{1},\cdots,x_{n})=(x_{1},\cdots,x_{n},0,\cdots,0),\\
\iota_{n+1}(x_{1},\cdots,x_{n+1})=(x_{1},\cdots,x_{n},0,\cdots,0,x_{n+1}).
\end{array}
\]

We first prove that $(i_{+1})_{\ast}$ is surjective for $0\leq i\leq n$.
Indeed, let $[a]\in\pi_{i}(\Omega(n+1,p+1))$\ be represented by $a:(S^{i}%
,\mathbf{e}_{1})\rightarrow(\Omega(n+1,p+1),\iota_{n+1})$. Then by Proposition
5.3 we have a homotopy $\varphi_{\lambda}:(S^{i},\mathbf{e}_{1})\rightarrow
(\Omega(n+1,p+1),\iota_{n+1})$\ such that $\varphi_{1}(S^{i})\subset
i_{+1}(\Omega(n,p))$.

Next let $[b]\in\pi_{i}(\Omega(n,p))$\ be represented by $b:(S^{i}%
,\mathbf{e}_{1})\rightarrow(\Omega(n,p),\iota_{n})$ such that $(i_{+1})_{\ast
}([b])=0$. Then we have a homotopy $\widetilde{\varphi}:S^{i}\times
I\rightarrow(\Omega(n+1,p+1),\iota_{n+1})$\ such that $\widetilde{\varphi
}|S^{i}\times0=i_{+1}\circ b$\ under the identification $S^{i}=S^{i}\times0$
and $\widetilde{\varphi}(\mathbf{e}_{1}\times I\cup S^{i}\times1)=\iota_{n+1}%
$. It follows from Proposition 5.3\ that if $i<n$, then there exists a
homotopy $\Phi_{\lambda}:(S^{i}\times I,S^{i}\times0\cup\mathbf{e}_{1}\times
I\cup S^{i}\times1)$\ relative to $S^{i}\times0\cup\mathbf{e}_{1}\times I\cup
S^{i}\times1$\ such that $\Phi_{1}(S^{i}\times I)\subset i_{+1}(\Omega(n,p))$.
This prove the injectivity of $(i_{+1})_{\ast}$.
\end{proof}

\begin{proposition}
Let (i) $n<p$, or (ii) $n=p>1$ and $I=(1,0)$. Then $(i_{+1})_{\ast}:\pi
_{i}(\Omega^{I}(n,p))\rightarrow\pi_{i}(\Omega^{I}(n+1,p+1))$ is an
isomorphism for $0\leq i<n$ and an epimorphism for $i=n$.
\end{proposition}

\begin{proof}
The case (i) follows from Lemma 5.1 and Proposition 5.4.

The case (ii) is proved as follows. Let $d:\mathbb{R}^{n+1}\rightarrow
\mathbb{R}^{n+1}$ be the diffeomorphism defined by $d(x_{1},\cdots
,x_{n+1})=(x_{n+1},x_{1},\cdots,x_{n})$ and let $d^{\Omega}:\Omega
^{(1,0)}(n+1,n+1)\rightarrow\Omega^{(1,0)}(n+1,n+1)$ be the diffeomorphism,
which maps a $k$-jet $j_{0}^{k}f$\ to $j_{0}^{k}(d\circ f\circ d^{-1})$.\ Let
$j:SO(n+1)\rightarrow SO(n+2)$ be the map defined by $j(A)=(1)\dotplus A$ for
$A\in SO(n+1)$. Let us recall the topological embedding $i_{m}^{SO}%
:SO(m+1)\rightarrow\Omega^{(1,0)}(m,m)$ defined in \cite[Proposition
2.4]{FoldSurg}, which is equivariant with respect to the action of
$SO(m)\times SO(m)$ and that $i_{m}^{SO}(SO(m+1))$ is a deformation retract of
$\Omega^{(1,0)}(m,m)$. Then we have the following commutative diagram. The
proof of commutativity is a rather wearisome calculation using the definition
of $i_{m}^{SO}$, and so it is left to the reader.%
\[%
\begin{array}
[c]{ccccc}%
SO(n+1) &  & \overset{j}{\longrightarrow} &  & SO(n+2)\\
{\scriptstyle{i_{n+1}^{SO}}}\downarrow\text{ \ \ } &  &  &  & \text{
\ \ \ \ }\downarrow{\scriptstyle{i_{n+2}^{SO}}}\\
\Omega^{(1,0)}(n,n) & \longrightarrow & \Omega^{(1,0)}(n+1,n+1) &
\overset{d^{\Omega}}{\longrightarrow} & \Omega^{(1,0)}(n+1,n+1),
\end{array}
\]
Then the assertion of the corollary follows from the corresponding assertion
for the map $j.$
\end{proof}

Recall that $\boldsymbol{\Omega}_{n}=\Omega(\gamma_{G_{n}}^{n},TP)$\ and
$\boldsymbol{\Omega}_{n+1}=\Omega(\gamma_{G_{n+1}}^{n+1},TP\oplus
\varepsilon_{P}^{1})$. Then Theorem 1.2 follows from Lemma 5.1, Proposition
5.5 and the following theorem.

\begin{theorem}
Let $\Omega(n,p)$\ and $\Omega(n+1,p+1)$\ denote the open subspaces, which
satisfy (\textbf{C}1), (\textbf{C}2) and (\textbf{C}3) as in Proposition 5.4
when $n<p$, or $\Omega^{(1,0)}(n,n)$ and $\Omega^{(1,0)}(n+1,n+1)$ when $n=p$
respectively. Let $G_{n}$ and $G_{n+1}$ refers to $G_{n,\ell}$ and
$G_{n+1,\ell}$ in the nonoriented case (resp. $\widetilde{G}_{n,\ell}$ and
$\widetilde{G}_{n+1,\ell}$ in the oriented case). Then the homomorphism%
\[
T(\mathbf{b(}\widehat{\gamma})^{(\Omega_{n},\Omega_{n+1})})_{\ast}:\pi
_{i+\ell}\left(  T(\widehat{\gamma}_{\boldsymbol{\Omega}_{n}}^{\ell})\right)
\longrightarrow\pi_{i+\ell}\left(  T(\widehat{\gamma}_{\boldsymbol{\Omega
}_{n+1}}^{\ell})\right)
\]
is an isomorphism for $0\leq i<n$ and an epimorphism for $i=n$.
\end{theorem}

\begin{proof}
We first show that $(i^{G})_{\ast}:\pi_{i}(G_{n})\rightarrow\pi_{i}(G_{n+1})$
is an isomorphism for $0\leq i<n$ and an epimorphism for $i=n$. We give a
proof only in the nonoriented case, since the proof in the oriented case is
analogous. Let us consider the canonical maps%
\begin{align*}
q_{1}  &  :O(n+\ell+1)/O(\ell)\times O(n)\times E_{1}\longrightarrow
G_{n+1,\ell},\\
q_{2}  &  :O(n+\ell+1)/O(\ell)\times O(n)\times E_{1}\longrightarrow
O(n+\ell+1)/O(n+\ell)\times E_{1}=S^{n+\ell},\\
j  &  :G_{n,\ell}\rightarrow O(n+\ell+1)/O(\ell)\times O(n)\times E_{1}.
\end{align*}
It is obvious that $q_{1}$ and $q_{2}$ yield the structures of the fiber
bundles with the fibers $S^{n}$\ and $G_{n,\ell}$, which is included by
$j$,\ respectively. Since $j_{\ast}:\pi_{i}(G_{n,\ell})\rightarrow\pi
_{i}(O(n+\ell+1)/O(\ell)\times O(n)\times E_{1})$ is an isomorphism for $0\leq
i<n+\ell-1$ and since $(q_{1})_{\ast}:\pi_{i}(O(n+\ell+1)/O(\ell)\times
O(n)\times E_{1})\rightarrow\pi_{i}(G_{n+1,\ell})$ is an isomorphism for
$0\leq i<n$ and an epimorphism for $i=n$, the assertion follows.

Let us recall the fiber map $\mathbf{j}^{(\Omega_{n},\Omega_{n+1})}$\ covering
$i^{G}\times id_{P}$. We next prove that%
\begin{equation}
(\mathbf{j}^{(\Omega_{n},\Omega_{n+1})})_{\ast}:\pi_{i}(\boldsymbol{\Omega
}_{n})\rightarrow\pi_{i}(\boldsymbol{\Omega}_{n+1})
\end{equation}

\noindent is an isomorphism for $0\leq i<n$ and an epimorphism for $i=n$. Let
us consider the diagram which is induced from the homomorphisms $(\mathbf{j}%
^{(\Omega_{n},\Omega_{n+1})})_{\ast}$\ of the exact sequence of\ the homotopy
groups for the fiber bundle $\boldsymbol{\Omega}_{n}$\ over $G_{n,\ell}\times
P$\ to the exact sequence of the fiber bundle $\boldsymbol{\Omega}_{n+1}%
$\ over $G_{n+1,\ell}\times P$. Then the second assertion\ about (5.3) follows
from \cite[Lemma 3.2]{CaEilen}, the first assertion about $(i^{G})_{\ast}$ and
Proposition 5.4. Then it follows from \cite[Section 5, 9 Theorem]{SpaAT} that%
\begin{equation}
(\mathbf{j}^{(\Omega_{n},\Omega_{n+1})})_{\ast}:H_{i}(\boldsymbol{\Omega}%
_{n})\rightarrow H_{i}(\boldsymbol{\Omega}_{n+1})
\end{equation}

\noindent is an isomorphism for $0\leq i<n$ and an epimorphism for $i=n$. By
virtue of the Thom Isomorphism Theorem, we have that%
\[
T(\mathbf{b(}\widehat{\gamma})^{(\Omega_{n},\Omega_{n+1})})_{\ast}:H_{i+\ell
}\left(  T(\widehat{\gamma}_{\boldsymbol{\Omega}_{n}}^{\ell})\right)
\longrightarrow H_{i+\ell}\left(  T(\widehat{\gamma}_{\boldsymbol{\Omega
}_{n+1}}^{\ell})\right)
\]
is an isomorphism for $0\leq i<n$ and an epimorphism for $i=n$. Since both of
Thom spaces are simply connected, the assertion of theorem follows from the
Hurewicz Isomorphism Theorem.
\end{proof}

Let $\mathcal{K}^{k}$ denote the group of\ $k$-jets of contact
transformations, which acts on $J^{k}(m,q)$, defined in \cite[(2.6)]{MathIII}.
The $\mathcal{K}^{k}$-orbit of $z\in J^{k}(m,q)$ is denoted by $\mathcal{K}%
^{k}z$. It is known \cite[Theorem 2.1]{MathIV} that $z^{\prime}$ lies in
$\mathcal{K}^{k}z$ if and only if $Q(z^{\prime})\approx Q(z)$. Hence, given a
$\mathbb{R}$-algebra $Q$ which has a $k$-jet $z$ with $Q\approx Q(z)$, we
write $\Sigma_{Q}=\mathcal{K}^{k}z$. In general, $Q(z)$ is written as
$\mathbb{R}[x_{1},\cdots,x_{m}]/(f_{1},\cdots,f_{q})+\mathfrak{m}_{m}^{k+1},$
and hence let $b$ denote the minimal number of generators for the ideal
$(f_{1},\cdots,f_{q})$ modulo $\mathfrak{m}_{m}^{k+1}$. If we define the
integer $\iota(Q)=m-b$, then it is easy to see that $\iota(Q)$ is an invariant
of the isomorphism class of $Q$. Du Plessis\cite{duPContIn} has proved the
following theorem (see also several examples which are concerned with simple singularities).

\begin{theorem}
[17]Let $\Omega(m,q)$ be an open set of $J^{k}(m,q)$, which is invariant with
respect to the action of $\mathcal{K}^{k}$. Assume that for each
$\mathcal{K}^{k}$-orbit $\Sigma_{Q}\subset\Omega(m,q)$, there exists a
$\mathcal{K}^{k}$-orbit $\Sigma_{Q^{\prime}}\subset\Omega(m,q)$ for another
$\mathbb{R}$-algebra $Q^{\prime}$ such that $\Sigma_{Q^{\prime}}%
\subset\mathrm{Closure}(\Sigma_{Q})$ and $-\iota(Q^{\prime})<q-m$. Then
$\Omega(m,q)$ is extensible.
\end{theorem}

\begin{lemma}
Assume that $\Omega(n,p)$\ and $\Omega(n+1,p+1)$\ satisfy (\textbf{C}1),
(\textbf{C}2) and (\textbf{C}3) and that $\Omega(n+1,p+1)\cap\Sigma
^{n+1}(n+1,p+1)=\emptyset$. Assume that $\Omega(n,p)$\ satisfies the
assumption of Theorem 5.7 for $(m,q)=(n,p)$. Then $\Omega(n+1,p+1)$ also
satisfies the assumption of Theorem 5.7 for $(m,q)=(n+1,p+1)$. Namely,
$\Omega(n,p)$\ and $\Omega(n+1,p+1)$\ are extensible.

In particular, the pair $(\Omega(n,p),\Omega(n+1,p+1))$\ is admissible to the h-Principle.
\end{lemma}

\begin{proof}
Take a $\mathcal{K}^{k}$-orbit $\Sigma_{Q}\subset\Omega(n+1,p+1)$. Then there
exists $k$-jets $z_{1}\in\Sigma_{Q}$ and $z_{2}\in\Omega(n,p)$ with $Q\approx
Q(z_{1})\approx Q(z_{2})$ by (\textbf{C}3). Hence, $\mathcal{K}^{k}%
z_{2}\subset\Omega(n,p)$. By assumption, there exists a $\mathcal{K}^{k}%
$-orbit $\Sigma_{Q^{\prime}}$ for another $\mathbb{R}$-algebra $Q^{\prime}$
such that $\Sigma_{Q^{\prime}}\subset$Closure$(\mathcal{K}^{k}z_{2})$ and
$-\iota(Q^{\prime})<p-n$. By (\textbf{C}1), $i_{+1}(\mathcal{K}^{k}%
z_{2})\subset\Omega(n+1,p+1)$ and $i_{+1}(\Sigma_{Q^{\prime}})\subset
i_{+1}(\mathrm{Closure}(\mathcal{K}^{k}z_{2}))$. Since $Q(z_{2})\approx Q$, we
have $\mathcal{K}^{k}(i_{+1}(z_{2}))=\Sigma_{Q}$, and hence we have%
\[
i_{+1}(\mathrm{Closure}(\mathcal{K}^{k}z_{2}))\subset\mathrm{Closure}%
(\mathcal{K}^{k}(i_{+1}(z_{2})))=\mathrm{Closure}(\Sigma_{Q}).
\]
Since the $\mathbb{R}$-algebras associated to $i_{+1}(\Sigma_{Q^{\prime}})$
are all isomorphic to $Q^{\prime}$, the $\mathcal{K}^{k}$-orbit of
$i_{+1}(\Sigma_{Q^{\prime}})$ in $J^{k}(n+1,p+1)$ is contained in
Closure$(\Sigma_{Q})$. This is what we\ want to prove.
\end{proof}

If we note codim$\Sigma^{n+1}(n+1,p+1)=(n+1)(p+1)$, we have the following
corollary of Theorems 4.2, 5.6 and Lemma 5.8.

\begin{corollary}
Let $n<p$. Assume that $\Omega(n,p)$\ and $\Omega(n+1,p+1)$\ satisfy
(\textbf{C}1), (\textbf{C}2) and (\textbf{C}3). Assume that $\Omega
(n,p)$\ satisfies the assumption of Theorem 5.7 for $(m,q)=(n,p)$. Then the
homomorphisms%
\begin{align*}
\omega_{\mathfrak{N}}^{(\Omega_{n},\Omega_{n+1})}  &  :\mathfrak{N{Cob}}%
_{n,P}^{(\Omega_{n},\Omega_{n+1})}\longrightarrow\pi_{n+\ell}\left(
T(\widehat{\gamma}_{\Omega(\gamma_{G_{n+1,\ell}}^{n+1},TP\oplus\varepsilon
_{P}^{1})}^{\ell})\right) \\
\omega_{\mathfrak{O}}^{(\Omega_{n},\Omega_{n+1})}  &  :\mathfrak{{OCob}}%
_{n,P}^{(\Omega_{n},\Omega_{n+1})}\longrightarrow\pi_{n+\ell}\left(
T(\widehat{\gamma}_{\Omega(\gamma_{\widetilde{G}_{n+1,\ell}}^{n+1}%
,TP\oplus\varepsilon_{P}^{1})}^{\ell})\right)
\end{align*}
are isomorphisms.
\end{corollary}

\section{Cobordisms of fold-maps}

In this section we study $\mathfrak{OCob}_{n,P}^{(\Omega^{(1,0)}%
,\Omega^{(1,0)})}$ in the equidimension $n=p$, where $G_{n+1}$ refers to
$\widetilde{G}_{n+1,\ell}$.\ For this we study $\pi_{n+\ell}(\widehat{\gamma
}_{\boldsymbol{\Omega}_{n+1}^{(1,0)}}^{\ell})$ in place by Theorems 1.1 and
1.2 and prove Theorem 1.3 . Take a Riemannian metric on $P$.

Let $SO_{n+2}(\gamma_{G_{n+1}}^{n+1}\oplus\varepsilon_{G_{n+1}}^{1}%
,TP\oplus\varepsilon_{P}^{2})$\ denote the total space of the open subbundle
of $\mathrm{Hom}(\gamma_{G_{n+1}}^{n+1}\oplus\varepsilon_{G_{n+1}}%
^{1},TP\oplus\varepsilon_{P}^{2})$, which is associated to $SO(n+2)$. Namely,
it consists of all isomorphisms which preserve orientations and norms of
vectors of fibers. Since the topological embedding $i_{n+1}^{SO}%
:SO(n+2)\rightarrow\Omega^{(1,0)}(n+1,n+1)$\ in \cite[Proposition
2.4]{FoldSurg} is equivariant with respect to the action of $SO(n+1)\times
SO(n+1)$ and since $i_{n+1}^{SO}(SO(n+2))$ is a deformation retract of
$\Omega^{(1,0)}(n+1,n+1)$, there exists the homotopy equivalent fiber map%
\[
\mathfrak{i}_{SO}:SO_{n+2}(\gamma_{G_{n+1}}^{n+1}\oplus\varepsilon_{G_{n+1}%
}^{1},TP\oplus\varepsilon_{P}^{2})\longrightarrow\boldsymbol{\Omega}%
_{n+1}^{(1,0)}=\Omega^{(1,0)}(\gamma_{G_{n+1}}^{n+1},TP\oplus\varepsilon
_{P}^{1})
\]
over $G_{n+1}\times P$ associated to $i_{n+1}^{SO}$. The image is denoted by
$\mathbf{SO}_{n+2}$. Furthermore, $\mathfrak{i}_{SO}$ induces the bundle map
\[
\mathfrak{B}:\widehat{\gamma}_{\boldsymbol{\Omega}_{n+1}^{(1,0)}}^{\ell
}|_{\mathbf{SO}_{n+2}}\longrightarrow\widehat{\gamma}_{\boldsymbol{\Omega
}_{n+1}^{(1,0)}}^{\ell}%
\]
covering the inclusion $\mathbf{SO}_{n+2}\rightarrow\boldsymbol{\Omega}%
_{n+1}^{(1,0)}$.

\begin{proposition}
(i) The fiber map $\mathfrak{i}_{SO}$ is a homotopy equivalence.

(ii) The Thom map%
\[
T(\mathfrak{B)}:T\left(  \widehat{\gamma}_{\boldsymbol{\Omega}_{n+1}^{(1,0)}%
}^{\ell}|_{\mathbf{SO}_{n+2}}\right)  \longrightarrow T\left(  \widehat
{\gamma}_{\boldsymbol{\Omega}_{n+1}^{(1,0)}}^{\ell}\right)
\]
is a homotopy equivalence.
\end{proposition}

Let $SO(TP\oplus\varepsilon_{P}^{1})$\ denote the total space of the principal
bundle over $P$ associated to $TP\oplus\varepsilon_{P}^{1}$, whose fiber is
$SO(n+1)$. Let $(SO(\ell)\times E_{n+1})\backslash SO(n+\ell+1)$ be the
Stiefel manifold. Consider the natural actions of $SO(n+1)$ on $SO(TP\oplus
\varepsilon_{P}^{1})$\ from the right-hand side and on $SO(n+2)$\ through
$SO(n+1)\times E_{1}$\ from the left-hand side, and the natural actions of
$SO(n+1)$ on $SO(n+2)$\ through $SO(n+1)\times E_{1}$\ from the right-hand
side and on $(SO(\ell)\times E_{n+1})\backslash SO(n+\ell+1)$ from the
left-hand side respectively. Then we can express as%
\begin{align*}
\mathbf{SO}_{n+2}  &  =SO(TP\oplus\varepsilon_{P}^{1})\underset{SO(n+1)}%
{\times}SO(n+2)\\
&  \text{ \ \ \ \ \ \ \ \ }\underset{SO(n+1)}{\times}((SO(\ell)\times
E_{n+1})\backslash SO(n+\ell+1)).
\end{align*}
Identify the quotient space $SO(TP\oplus\varepsilon_{P}^{1})/SO(n+1)$\ with
$P$. Then we have the projection $pr_{P}^{SO}:\mathbf{SO}_{n+2}\rightarrow
P$\ by forgetting the component $SO(n+2)\underset{SO(n+1)}{\times}%
((SO(\ell)\times E_{n+1})\backslash SO(n+\ell+1))$, denoted by $\mathfrak{f}$,
which is the canonical fiber of $pr_{P}^{SO}$. Let $pr_{S^{n+1}}%
^{\mathfrak{{f}}}:\mathfrak{f}\rightarrow SO(n+2)/SO(n+1)=S^{n+1}$ be the
projection forgetting the component $(SO(\ell)\times E_{n+1})\backslash
SO(n+\ell+1)$. Since the last space is $(\ell-1)$-connected,%
\begin{equation}
(pr_{S^{n+1}}^{\mathfrak{f}})_{\ast}:\pi_{i}(\mathfrak{f})\rightarrow\pi
_{i}(S^{n+1})
\end{equation}
is an isomorphism for $i<\ell$. Hence, we have the following lemma.

\begin{lemma}
The homomorphism $(pr_{P}^{SO})_{\ast}:\pi_{i}(\mathbf{SO}_{n+2}%
)\rightarrow\pi_{i}(P)$ is an isomorphism for $i<n+1$ and is an epimorphism
for $i=n+1$.
\end{lemma}

Let us define a bundle map%
\[
\mathbf{b}_{\mathcal{T}}:(\pi_{G_{n+1}}^{k})^{\ast}(\gamma_{G_{n+1}}%
^{n+1})|_{\mathbf{SO}_{n+2}}\longrightarrow TP\oplus\varepsilon_{P}^{1}%
\]
as follows. We denote an element of $\mathbf{SO}_{n+2}$ and $(\pi_{G_{n+1}%
}^{k})^{\ast}(\gamma_{G_{n+1}}^{n+1})|_{\mathbf{SO}_{n+2}}$ by $[a,y,h]$\ and
by $[a,y,h,\mathbf{v]}$ respectively, where $a\in G_{n+1}$, $y\in P$,
$\mathbf{v\in}a$\ and $h:a\rightarrow T_{y}P\oplus\mathbb{R}$\ is an
isomorphism preserving orientations and norms. Then we set $\mathbf{b}%
_{\mathcal{T}}([a,y,h,\mathbf{v])}=h(\mathbf{v)}$.\ Let us consider the
trivialization $(t_{P}\oplus id_{P}\times\mathbb{R})\circ(id_{TP}%
\oplus\mathbf{k}_{P}^{\backsim}):TP\oplus\varepsilon_{P}^{1}\oplus\nu
_{P}^{\ell}\rightarrow\varepsilon_{P}^{n+\ell+1}$. By \cite[Proposition
3.3]{FoldSurg}, there exists a bundle map%
\begin{equation}
\mathbf{b}_{\nu}:\widehat{\gamma}_{\boldsymbol{\Omega}_{n+1}^{(1,0)}}^{\ell
}|_{\mathbf{SO}_{n+2}}\rightarrow\nu_{P}^{\ell}%
\end{equation}
such that $(t_{P}\oplus id_{P}\times\mathbb{R})\circ(id_{TP}\oplus
\mathbf{k}_{P}^{\backsim})\circ(\mathbf{b}_{\mathcal{T}}\oplus\mathbf{b}_{\nu
})\circ t_{\boldsymbol{\Omega}_{n+1}^{(1,0)}}|_{\mathbf{SO}_{n+2}}$\ is
homotopic to $pr_{P}^{SO}\times id_{\mathbb{R}^{n+\ell+1}}$. Then we have the
following lemma.

\begin{lemma}
The homomorphism induced from the Thom map of $\mathbf{b}_{\nu}$%
\[
T(\mathbf{b}_{\nu})_{\ast}:\pi_{i+\ell}\left(  T(\widehat{\gamma
}_{\boldsymbol{\Omega}_{n+1}^{(1,0)}}^{\ell}|_{\mathbf{SO}_{n+2}})\right)
\longrightarrow\pi_{i+\ell}\left(  T(\nu_{P}^{\ell})\right)
\]
is an isomorphism for $i<n+1$ and is an epimorphism for $i=n+1$.
\end{lemma}

\begin{proof}
By Lemma 6.2, (6.2), \cite[Section 5, 9 Theorem]{SpaAT} and the Thom
Isomorphism Theorem,%
\[
T(\mathbf{b}_{\nu})_{\ast}:H_{i+\ell}\left(  T(\widehat{\gamma}%
_{\boldsymbol{\Omega}_{n+1}^{(1,0)}}^{\ell}|_{\mathbf{SO}_{n+2}})\right)
\longrightarrow H_{i+\ell}\left(  T(\nu_{P}^{\ell})\right)
\]
is an isomorphism for $i<n+1$ and is an epimorphism for $i=n+1$. Since
$T((\pi_{G_{n+1}}^{k})^{\ast}(\widehat{\gamma}_{G_{n+1}}^{\ell})|_{\mathbf{SO}%
_{n+2}})$ and $T(\nu_{P}^{\ell})$ are simply connected, it follows from the
Hurewicz Isomorphism Theorem that%
\[
T(\mathbf{b}_{\nu})_{\ast}:\pi_{i+\ell}\left(  T(\widehat{\gamma
}_{\boldsymbol{\Omega}_{n+1}^{(1,0)}}^{\ell}|_{\mathbf{SO}_{n+2}})\right)
\longrightarrow\pi_{i+\ell}\left(  T(\nu_{P}^{\ell})\right)
\]
is an isomorphism for $i<n+1$ and is an epimorphism for $i=n+1$.
\end{proof}

According to \cite{SpaDual}, let $\{S^{n+\ell},T(\nu_{P}^{\ell})\}$ denote the
set of S-homotopy classes of S-maps $S^{i}\wedge S^{n+\ell}\rightarrow
S^{i}\wedge T(\nu_{P}^{\ell})$ $({i}\geq0)$. An element of $\{S^{n+\ell}%
,T(\nu_{P}^{\ell})\}$ represented by a map $\alpha:S^{n+\ell}\rightarrow
S^{i}\wedge T(\nu_{P}^{\ell})$ is written as $\{\alpha\}$. Since $\ell\gg n$,
$\{S^{n+\ell},T(\nu_{P}^{\ell})\}$ is isomorphic to $\pi_{n+\ell}\left(
T(\nu_{P}^{\ell})\right)  $. It has been proved in \cite[Lemma 2]{MiSpa} that
$T(\nu_{P}^{\ell})$ is the S-dual space of $P^{0}=P\cup\ast_{P}$, where
$\ast_{P}$ is the base point. Namely, we have the isomorphism $\{S^{n+\ell
},T(\nu_{P}^{\ell})\}\approx\{S^{\ell}P^{0},S^{\ell}\}$.

Let $\mathfrak{OCob}_{n,P}^{(\Omega^{(1,0)},\Omega^{(1,0)})}(d)$ denote the
subset of $\mathfrak{OCob}_{n,P}^{(\Omega^{(1,0)},\Omega^{(1,0)})}$ which
consists of all fold-maps into $P$ of degree $d$. In particular,
$\mathfrak{OCob}_{n,P}^{(\Omega^{(1,0)},\Omega^{(1,0)})}(0)$ is a submodule of
$\mathfrak{OCob}_{n,P}^{(\Omega^{(1,0)},\Omega^{(1,0)})}$. Let $\{S^{n+\ell
},T(\nu_{P}^{\ell})\}_{d}$\ denote the subset of $\{S^{n+\ell},T(\nu_{P}%
^{\ell})\}$ which consists of all S-maps of degree $d$. By Theorem 1.1,
Proposition 6.1 and Lemma 6.3, $\mathfrak{OCob}_{n,P}^{(\Omega^{(1,0)}%
,\Omega^{(1,0)})}(d)$ is mapped bijectively onto $\{S^{n+\ell},T(\nu_{P}%
^{\ell})\}_{d}$.

We prove the following refined form of Theorem 1.3.

\begin{theorem}
Let $n=p\geq2$ and $P$ be a closed, connected, oriented $n$-dimensional
manifold. Then there exist a bijection $\omega_{d}:\mathfrak{OCob}%
_{n,P}^{(\Omega^{(1,0)},\Omega^{(1,0)})}(d)\rightarrow\lbrack P,F^{d}]$.
\end{theorem}

\begin{proof}
Let us define the map $c_{F}:\{S^{\ell}P^{0},S^{\ell}\}\rightarrow\lbrack
P,F]$. Let $\{\beta\}\in\{S^{\ell}P^{0},S^{\ell}\}$ be represented by
$\beta:S^{\ell}P^{0}\rightarrow S^{\ell}$. For a point $x\in P$ we define
$\beta(x):S^{\ell}=S^{0}\wedge S^{\ell}\rightarrow S^{\ell}$ by $(\beta
|\{\ast_{P}\cup x\}\wedge S^{\ell})\circ(\iota_{x}\wedge id_{S^{\ell}})$,
where $\iota_{x}:S^{0}\rightarrow\{\ast_{P}\cup x\}$ is the canonical
identification. Then we set $c_{F}(\{\beta\})(x)=\{\beta(x)\}$. It is easy to
check that $c_{F}$ is bijective. Furthermore, we have proved in \cite[Lemma
2.4]{FoldKyoto} that $c_{F}$ maps $\{S^{n+\ell},T(\nu_{P}^{\ell})\}_{d}$ to
the subset of $\{S^{\ell}P^{0},S^{\ell}\}$ which consists of all $\{\beta\}$
such that $\beta(x)$ is of degree $d$, namely $c_{F}(\{\beta\})(x)\in F_{d}$
for any $x\in P$. This proves the assertion.
\end{proof}

Let $\pi_{n}^{S}$ denote\ the $n$-th stable homotopy group of spheres
$\lim_{\ell\rightarrow\infty}\pi_{n+\ell}(S^{\ell})$. It follows from
\cite{AtThomCom} that $[S^{n},F^{0}]$ is canonically isomorphic to $\pi
_{n}^{S}$. So identifying $[S^{n},F^{0}]$ with $\pi_{n}^{S}$, we have the
following corollary.

\begin{corollary}
The map $\omega_{0}:\mathfrak{OCob}_{n,S^{n}}^{(\Omega^{(1,0)},\Omega
^{(1,0)})}(0)\rightarrow\pi_{n}^{S}$ is an isomorphism.
\end{corollary}

\section{Stable maps of spheres}

Let $\mathfrak{F}^{(\Omega^{(1,0)},\Omega^{I})}:\mathfrak{OCob}_{n,S^{n}%
}^{(\Omega^{(1,0)},\Omega^{(1,0)})}(0)\rightarrow\mathfrak{OCob}_{n,S^{n}%
}^{(\Omega^{(1,0)},\Omega^{I})}(0)$ denote the homomorphism which maps an
$\Omega^{(1,0)}$-cobordism class $[f]$ to the $\Omega^{I}$-cobordism class of
$f$. Let $I([f])$ denote the smallest symbol $I$ such that $\mathfrak{F}%
^{(\Omega^{(1,0)},\Omega^{I})}([f])$ is a null element. Then there exists an
$\Omega^{I([f])}$-regular cobordism $F:(V,\partial V)\rightarrow(S^{n}\times
I,S^{n}\times0)$ such that $\partial V=N$, the collar of $\partial V$ is
identified with $N\times\lbrack0,\varepsilon]$,and $F|N\times\lbrack
0,\varepsilon]=f\times id_{[0,\varepsilon]}$. In this section we show that the
singularities of symbol $I([f])$\ of $F$ detect the stable homotopy class
$\omega_{0}([f])\in\pi_{n}^{S}$ in low dimensions. We have to prepare some
machinery for this purpose, although the dimensions are low.

Let $\mathbf{D}$ and $\mathbf{P}$ denote the total tangent bundle defined on
$J^{\infty}(V,Y)$ and $(\pi_{Y}^{\infty})^{\ast}(TY)$ respectively. Let us
recall the fundamental property of $\mathbf{D}$\ over $J^{\infty}(V,Y)$. Let
$f:(V,x)\rightarrow(Y,y)$ be a map defined on a neighborhood\ $U_{x}$\ of
$x$\ with coordinates $(x_{1},\cdots,x_{n+1})$ and $\digamma$ be a smooth
function in the sense of \cite[Definition 1.4]{Board} defined on a
neighborhood of $j_{x}^{\infty}f$. We have the local vector fields $D_{i}%
$\ defined around $z$ with the property%
\[
D_{i}\digamma\circ j^{\infty}f=\frac{\partial}{\partial x_{i}}(\digamma\circ
j^{\infty}f)\text{ \ \ }(1\leq i\leq n+1),
\]
which span $\mathbf{D}$. It follows that $d(j^{\infty}f)(\partial/\partial
x_{i})(\digamma)=D_{i}\digamma(j^{\infty}f)$, where $d(j^{\infty
}f):TV\rightarrow T(J^{\infty}(V,Y))$ around $x$. This implies $d(j^{\infty
}f)(\partial/\partial x_{i})=D_{i}.$ Hence, we have $\mathbf{D\cong(}\pi
_{V}^{\infty})^{\ast}(TV)$. There have been defined the homomorphism
$\mathbf{d}_{1}:\mathbf{D}\rightarrow\mathbf{P}$ over $J^{\infty}(V,Y)$. If
$z=j_{x}^{\infty}f$, then $\mathbf{d}_{1,z}(D_{i})=(z,d_{x}f(\partial/\partial
x_{i}))$. The manifold $\Sigma^{1}(V,Y)$ is defined to be the submanifold of
$J^{\infty}(V,Y)$ which consists of all jets $z$ such that the kernel rank of
$\mathbf{d}_{1,z}$ is $1$. Since $\mathbf{d}_{1}|_{\Sigma^{1}(V,Y)}$ is of
constant rank $n$, we set $\mathbf{K}_{1}=$Ker$(\mathbf{d}_{1})$ and
$\mathbf{P}_{1}=$Cok$(\mathbf{d}_{1})$, which are vector bundles over
$\Sigma^{1}(V,Y)$. Let $\mathbf{1}_{r}$ denote $(\overbrace{1,\cdots,1}^{r})$.
The Boardman manifold $\Sigma^{\mathbf{1}_{r}}(V,Y)$ ($r\geq1$) has the
following properties (\cite{Board}).

(7-i) There exists the $(r+1)$-th intrinsic derivative%
\[
\mathbf{d}_{r+1}:T(\Sigma^{\mathbf{1}_{r-1}}(V,Y))|_{\Sigma^{\mathbf{1}_{r}%
}(V,Y)}\longrightarrow\mathrm{Hom}(S^{r}\mathbf{K}_{1},\mathbf{P}%
_{1})|_{\Sigma^{\mathbf{1}_{r}}(V,Y)}\longrightarrow\mathbf{0},
\]
so that $\mathrm{Ker}(\mathbf{d}_{r+1})=T(\Sigma^{\mathbf{1}_{r}}(V,Y))$.
Namely, $\mathbf{d}_{k+1}$ induces the isomorphism of the normal bundle
$(T(\Sigma^{\mathbf{1}_{r-1}}(V,Y))|_{\Sigma^{\mathbf{1}_{r}}(V,Y)}%
)/T(\Sigma^{\mathbf{1}_{r}}(V,Y))$ of $\Sigma^{\mathbf{1}_{r}}(V,Y)$ in
$\Sigma^{\mathbf{1}_{r-1}}(V,Y)$ onto $\mathrm{Hom}(S^{r}\mathbf{K}%
,\mathbf{P}_{1})|_{\Sigma^{\mathbf{1}_{r}}(V,Y)}$.

(7-ii) $\Sigma^{\mathbf{1}_{r+1}}(V,Y)$ is defined to be the submanifold of
$\Sigma^{\mathbf{1}_{r}}(V,Y)$ which consists of all jets $z$ such that
$\mathbf{d}_{r+1,z}|\mathbf{K}_{1,z}$ vanishes.

(7-iii) The $(r+2)$-th intrinsic derivative $\mathbf{d}_{r+2}$ is defined to
be the intrinsic derivative%
\[
d(\mathbf{d}_{r+1}|\mathbf{K}_{r}):T(\Sigma^{\mathbf{1}_{r}}(V,Y))|_{\Sigma
^{\mathbf{1}_{r+1}}(V,Y)}\rightarrow\text{Hom}(\mathbf{K}_{1},\mathbf{P}%
_{1})|_{\Sigma^{\mathbf{1}_{r+1}}(V,Y)}.
\]

(7-iv) The submanifold $\Sigma^{\mathbf{1}_{r}}(V,Y)$ is actually defined so
that it coincides with the inverse image of the submanifold $\widetilde
{\Sigma}^{\mathbf{1}_{r}}(V,Y)$ in $J^{r}(V,Y)$ by $\pi_{r}^{\infty}$. The
codimension of $\Sigma^{\mathbf{1}_{r}}(V,Y)$ in $J^{\infty}(V,Y)$ is $r$.

\begin{theorem}
Let $V$ be an oriented $(n+1)$-manifold with $\partial V$, which may be empty,
$Y$ be an oriented $(n+1)$-manifold and let $C$ be a closed subset of $V$. Let
$s$ be a section of $\Gamma_{\Omega^{1}}(V,Y)$ which has a fold-map $g$
defined on a neighborhood of $C$ into $Y$, where $j^{\infty}g=s$. Then there
exists an $\Omega^{(1,1,0)}$-regular map $f:V\rightarrow Y$ and a homotopy
$s_{\lambda}\in\Gamma_{\Omega^{1}}(V,Y)$\ relative to a neighborhood of $C$
such that $s_{0}=s$ and $s_{1}=j^{\infty}f.$
\end{theorem}

\begin{proof}
In the proof we use the notation introduced in \cite{Board}. By (2.2) we
always identify $J^{r}(V,Y)$ with $J^{r}(TV,TY)$, where $r$ may be $\infty$.
We may assume that $s$ is transverse to $\Sigma^{1}(V,Y)$ and we set
$S^{\mathbf{1}_{r}}(s)=s^{-1}(\Sigma^{\mathbf{1}_{r}}(V,Y))$. It follows that
$(\pi_{r}^{\infty}\circ s)(V\setminus(S^{\mathbf{1}_{r}}(s)))\subset
\Omega^{\mathbf{1}_{r-1},0}(V\setminus S^{\mathbf{1}_{r}}(s),Y)$.

We find a section $\mathfrak{s}$ of $\Omega^{(1,1,0)}(V,Y)$ such that $\pi
_{2}^{\infty}\circ\mathfrak{s}=\pi_{2}^{\infty}\circ s$. We set $(s|S^{1}%
(s))^{\ast}\mathbf{K}_{1}=K_{1}$ and $(s|S^{1}(s))^{\ast}\mathbf{P}_{1}=P_{1}%
$. Since $V$ and $Y$ are oriented and since $K_{1}$\ and $P_{1}$\ are line
bundles, we have that $K_{1}$\ and $P_{1}$\ are isomorphic. In particular, we
have the isomorphism $K_{1}|_{S^{\mathbf{1}_{2}}(s)}\rightarrow P_{1}%
|_{S^{\mathbf{1}_{2}}(s)}$. Consider the homomorphism%
\[
\mathbf{r}^{3}:\mathrm{Hom}(S^{3}(TV),TY)|_{S^{\mathbf{1}_{2}}(s)}%
\longrightarrow\mathrm{Hom}(S^{3}K_{1},P_{1})|_{S^{\mathbf{1}_{2}}(s)}%
\]
which is induced from the inclusion $S^{3}K_{1}|_{S^{\mathbf{1}_{2}}%
(s)}\rightarrow S^{3}(TV),TY)|_{S^{\mathbf{1}_{2}}(s)}$\ and the projection
$TY|_{S^{\mathbf{1}_{2}}(s)}\rightarrow P_{1}|_{S^{\mathbf{1}_{2}}(s)}$. Since
$S^{3}K_{1}\approx K_{1}$, there exists the isomorphism $\iota^{3}:S^{3}%
K_{1}|_{S^{\mathbf{1}_{2}}(s)}\rightarrow P_{1}|_{S^{\mathbf{1}_{2}}(s)}$,
which induces the isomorphism $K_{1}|_{S^{\mathbf{1}_{2}}(s)}\rightarrow
\mathrm{Hom}(S^{2}K_{1},P_{1})|_{S^{\mathbf{1}_{2}}(s)}.$ Since $S^{\mathbf{1}%
_{2}}(s)$\ is a closed submanifold of $V$\ such that $S^{\mathbf{1}_{2}%
}(s)\cap C=\emptyset$\ in $V$, there exists a homomorphism $\mathbf{h}%
^{3}:S^{3}(TV)|_{S^{\mathbf{1}_{2}}(s)}\rightarrow(\pi_{Y}^{\infty}\circ
s)^{\ast}(TY)|_{S^{\mathbf{1}_{2}}(s)}$\ such that $\mathbf{r}^{3}\circ
$\ $\mathbf{h}^{3}=\iota^{3}$. We extend $\mathbf{h}^{3}$ to the homomorphism
$\mathbf{H}^{3}:S^{3}(TV)\rightarrow(\pi_{Y}^{\infty}\circ s)^{\ast}(TY)$. If
$(\pi_{Y}^{\infty}\circ s)_{TY}:(\pi_{Y}^{\infty}\circ s)^{\ast}%
(TY)\rightarrow TY$\ denote the canonical bundle map covering $\pi_{Y}%
^{\infty}\circ s$, then we define the section $\mathfrak{s}:V\rightarrow
J^{\infty}(TV,TY)$ so that
\[
\pi_{3}^{\infty}\circ\mathfrak{s}(x)=\pi_{2}^{\infty}\circ s(x)\oplus(\pi
_{Y}^{\infty}\circ s)_{TY}\circ\mathbf{H}^{3}|_{x}%
\]
and that $\mathfrak{s}$\ is the composite of $\pi_{3}^{\infty}\circ
\mathfrak{s}$ and the canonical inclusion $J^{3}(TV,TY)\rightarrow J^{\infty
}(TV,TY)$.

We now show that $\mathfrak{s}(V)\subset\Omega^{(1,1,0)}(V,Y)$. In fact, it is
obvious that $\mathfrak{s}(V)\subset\Omega^{\mathbf{1}_{2}}(V,Y)$. It remains
to prove that if $x\in S^{\mathbf{1}_{2}}(s)$, then%
\begin{equation}
\mathbf{d}_{3,\mathfrak{s}(x)}:\mathbf{K}_{1,\mathfrak{s}(x)}\longrightarrow
\mathrm{Hom}(S^{2}\mathbf{K}_{1,\mathfrak{s}(x)},\mathbf{P}_{1,\mathfrak{s}%
(x)})
\end{equation}
is an isomorphism. In other words the homomorphism $S^{3}\mathbf{K}%
_{1,\mathfrak{s}(x)}\rightarrow\mathbf{P}_{1,\mathfrak{s}(x)}$ induced from
$\mathbf{d}_{3,\mathfrak{s}(x)}$ is an isomorphism. For any point $x\in
S^{\mathbf{1}_{2}}(s)$, let $y=\pi_{Y}^{\infty}\circ\mathfrak{s}(x)$, $U_{x}$
and $V_{y}$ be convex neighborhoods of $x$ and $y$ respectively. Let $t$\ and
$u\ $be the coordinates of $\exp_{V,x}(K_{1,x})$ and $\exp_{Y,y}((\pi
_{Y}^{\infty}\circ\mathfrak{s})_{TY}(P_{1,y}))$\ respectively, where $P_{1}$
is regarded as a line subbundle of $(\pi_{Y}^{\infty}\circ s)^{\ast
}(TY)|_{S^{\mathbf{1}_{2}}(s)}$ by virtue of the metric of $Y$. It follows
from (7.1), (2.2) and the definition of $\iota^{3}$ that%
\[%
\begin{array}
[c]{ll}%
(\bigcirc^{3}D_{t})u|_{\mathfrak{s}(x)}=\partial^{3}u/\partial t^{3}(x)\neq0 &
\text{for }x\in S^{\mathbf{1}_{2}}(s).
\end{array}
\]
Hence, we have that $\mathfrak{s}(S^{\mathbf{1}_{2}}(s))\subset\Sigma
^{(1,1,0)}(V,Y)$.

By \cite[Theorem 0.1]{hPrinBo} and \cite{G1} there exists an $\Omega
^{(1,1,0)}$-regular map $G:V\rightarrow Y$ such that $j^{\infty}G$ and
$\mathfrak{s}$ are homotopic relative to a neighborhood of $C$ as sections of
$\Omega^{1}(V,Y)$ over $V$. Here, we again note the remark which has been
given at the end of the proof of Theorem 1.1.
\end{proof}

Let $\ell$ be a natural number with $\ell\gg n.$\ Let $V$\ be an
$(n+1)$-manifold with $\partial V=N$, and let $\tau_{V}$\ be the stable $\ell
$-dimensional tangent bundle of $V$. Given a fold map $f:N\rightarrow S^{n}$
of degree $0$, we have the bundle map $\mathcal{T}(f):TN\oplus\varepsilon
_{N}^{1}\rightarrow\varepsilon_{S^{n}}^{n+1}$\ by \cite[Corollary 2]%
{FoldSurg}. Let us consider the obstruction for $\mathcal{T}(f)$\ to be
extended to the trivialization of $\tau_{V}$. For this purpose, we have the
primary obstruction $\mathfrak{o}(\tau_{V},\mathcal{T}(f))$ defined in
$H^{i+1}(V,N;\pi_{i}(SO(\ell)))$ for some $i$. Let $\widehat{V}=V\cup_{N}CN$,
which is obtained by pasting $V$ and the cone of $N$. Let $\tau(\widehat
{V},\mathcal{T}(f))$\ be the $\ell$-dimensional vector bundle, which is
obtained by pasting $\tau_{V}$\ and $\varepsilon_{CN}^{\ell}$\ by using
$\mathcal{T}(f)$. We have the primary obstruction $\mathfrak{o}(\tau
(\widehat{V},\mathcal{T}(f)))\in H^{i+1}(\widehat{V};\pi_{i}(SO(\ell)))\approx
H^{i+1}(V,N;\pi_{i}(SO(\ell)))$\ for $\tau(\widehat{V},\mathcal{T}(f))$\ to be
trivial.\ It is not difficult to see that $\mathfrak{o}(\tau_{V}%
,\mathcal{T}(f))=\pm\mathfrak{o}(\tau(\widehat{V},\mathcal{T}(f)))$ under the isomorphism.

\begin{remark}
In this case we may take $\ell=n+2$ and consider the subbundle $SO(\tau
(\widehat{V},\mathcal{T}(f)),\varepsilon_{\widehat{V}}^{n+2})$ of
$\mathrm{Hom}(\tau(\widehat{V},\mathcal{T}(f)),\varepsilon_{\widehat{V}}%
^{n+2})$\ associated to $SO(n+2)$. Since $i_{n+2}^{SO}:SO(n+2)\rightarrow
\Omega^{(1,0)}(n+1,n+1)$\ is a homotopy equivalence, $\mathfrak{o}%
(\tau(\widehat{V},\mathcal{T}(f)))$\ coincides with the obstruction to find a
section of $\Omega^{(1,0)}(\tau(\widehat{V},\mathcal{T}(f)),\varepsilon
_{\widehat{V}}^{n+2})$, namely the Thom polynomial of the closure
$Cl(\Sigma^{(1,1)}(\tau(\widehat{V},\mathcal{T}(f)),\varepsilon_{\widehat{V}%
}^{n+2}))$\ (see, for example, \cite[Proposition 3.1]{EliTB}). This Thom
polynomial is equal to the second Stiefel-Whitney class $w_{2}(\tau
(\widehat{V},\mathcal{T}(f)))$ by \cite{Porteou}.
\end{remark}

If $n+1=4m$ and $\omega_{0}([f])$ lies in what is called the $J$-image
$J(\pi_{4m-1}(SO(\ell)))$ of order $j_{m}$ in \cite{JXIV}, then we can choose
an fold-map $f$ such that $N=S^{n}$ by \cite[Proposition 5.1]{FoldKyoto}. This
is also true for the case $n=1$. Furthermore, we can take $V$ to be a
parallelizable manifold. Hence, the above dimension $i$ is equal to $n$.\ 

We have the following lemma due to \cite[Lemma 2]{MiKerBer}. Let $a_{m}$
denote $2$\ for $m$ odd and $1$ for $m$ even.

\begin{lemma}
[32]Let $n+1=4m$. Let $V$ be a parallelizable manifold. Then $\mathfrak{o}%
(\tau(\widehat{V},\mathcal{T}(f)))$ is related to the $m$-th Pontrjagin class
$P_{m}(\tau(\widehat{V},\mathcal{T}(f)))$ by the identity $P_{m}(\tau
(\widehat{V},\mathcal{T}(f)))=\pm a_{m}(2m-1)!\mathfrak{o}(\tau(\widehat
{V},\mathcal{T}(f)))$.
\end{lemma}

We next see how $\mathfrak{o}(\tau_{V},\mathcal{T}(f))$ varies depending on
the choice of $V$\ and $f$ (the following argument is available for the case
$n=1$). Let two fold maps $f_{i}:S^{n}\rightarrow S^{n}$\ of degree $0$
($i=0,1$) are $\Omega^{(1,0)}$-cobordant by a cobordism $F:(W,\partial
W)\rightarrow(S^{n}\times I,S^{n}\times0\cup S^{n}\times1)$\ of degree $0$ as
in Definition 3.1. Assume that there exists a parallelizable $(n+1)$-manifold
$V_{i}$\ with $\partial V_{i}=S^{n}\times i$. Then we have the bundle maps
$\mathcal{T}(f_{i}):T(S^{n}\times i)\oplus\varepsilon_{S^{n}\times i}%
^{1}\rightarrow\varepsilon_{S^{n}}^{n+1}$ and $\mathcal{T}(F):TW\oplus
\varepsilon_{W}^{1}\rightarrow\varepsilon_{S^{n}\times I}^{n+2}$\ by
\cite[Corollary 2]{FoldSurg} such that $TW|_{S^{n}\times i}=T(S^{n}\times
i)\oplus\varepsilon_{S^{n}\times i}^{1}$ and the stabilizations of
$\mathcal{T}(F)|_{N_{i}}$ and $\mathcal{T}(f_{i})$\ are equal. Consider the
almost parallelizable manifold $\widehat{W}=V_{0}\cup_{S^{n}\times0}%
W\cup_{S^{n}\times1}(-V_{1})$, which is obtained by pasting $V_{0}$, $W$ and
$V_{1}$\ with orientation reversed.\ Let%
\[
\mathfrak{o}(\tau(\widehat{W}))\in H^{n+1}(\widehat{W};\pi_{n}(SO(\ell
)))\approx\pi_{n}(SO(\ell))
\]
be the unique primary obstruction\ for $\tau(\widehat{W})$\ to be trivial.
Then it is obvious that%
\[
\mathfrak{o}(\tau_{V_{0}},\mathcal{T}(f_{0}))-\mathfrak{o}(\tau_{V_{1}%
},\mathcal{T}(f_{1}))=\pm\mathfrak{o}(\tau(\widehat{W}))\text{ \ \ in }\pi
_{n}(SO(\ell)).
\]
Define the integer $\mathfrak{m}(n)$ for $n>1$\ to be the minimal nonnegative
number such that there exists an $(n+1)$-dimensional almost parallelizable
closed manifold\ $\widehat{W}^{\prime}$\ such that $\mathfrak{o}(\tau
(\widehat{W}^{\prime}))=\mathfrak{m}(n)$. We will see $\mathfrak{m}(1)=0$
later. We have the following theorem due to \cite[Theorems 1 and 2]{MiKerBer}.

\begin{theorem}
[32]Let $n+1=4m$. Then we have

(i) The Pontrjagin class $P_{m}(\widehat{W}^{\prime})$\ of an almost
parallelizable closed manifold $\widehat{W}^{\prime}$\ is divisible by $\pm
j_{m}a_{m}(2m-1)!$.

(ii) There exists an almost parallelizable closed manifold $\widehat{W}_{0}$
with $P_{m}(\widehat{W}_{0})=\pm j_{m}a_{m}(2m-1)!$.
\end{theorem}

Consequently, we have $\mathfrak{m}(n)=j_{m}$.

Let $\mathbf{P}^{n}$ denote the real projective space of dimension $n$. Let
$C_{+}\mathbf{P}^{n}$ and $C_{-}\mathbf{P}^{n}$ denote the cone of
$\mathbf{P}^{n}\times\lbrack0,1]/\mathbf{P}^{n}\times1$ and $\mathbf{P}%
^{n}\times\lbrack-1,0]/\mathbf{P}^{n}\times(-1)$\ respectively and let
$S\mathbf{P}^{n}$\ denote the suspension $C_{+}\mathbf{P}^{n}\cup
_{\mathbf{P}^{n}}C_{-}\mathbf{P}^{n}$. For the projection of the double
covering $p_{\mathbf{P}^{n}}:S^{n}\rightarrow\mathbf{P}^{n}$, let
$S(p_{\mathbf{P}^{n}}):S^{n+1}\rightarrow S\mathbf{P}^{n}$\ be its suspension.

\begin{lemma}
Let $n+1=4m$. Then we have

(i) $S(p_{\mathbf{P}^{n}})_{\ast}:\pi_{i}(S^{n+1})\rightarrow\pi
_{i}(S\mathbf{P}^{n})$\ is an isomorphism modulo $2$-torsion for $0\leq i\leq
n+1$.

(ii) $S(p_{\mathbf{P}^{n}})_{\ast}:\pi_{n+1}(S^{n+1})\rightarrow\pi
_{n+1}(S\mathbf{P}^{n})$ is injective.
\end{lemma}

\begin{proof}
Since $S\mathbf{P}^{n}$\ is simply connected and $H_{i}(S\mathbf{P}^{n})$ is a
$2$-torsion for $0<i<n+1$, the assertion (i) follows from \cite[Section 6, 21
Theorem]{SpaAT}.

Let $\mathfrak{q}:S\mathbf{P}^{n}\rightarrow S(\mathbf{P}^{n}/\mathbf{P}%
^{n-1})=S^{n+1}$ be the collapsing map. Then it is obvious that the degree of
$\mathfrak{q}\circ S(p_{\mathbf{P}^{n}})$ is equal to $2$. This shows (ii).
\end{proof}

Let us define the map $\mathfrak{g}:\mathbf{P}^{n}\rightarrow SO(n+1)$ as
follows. Let $I_{-}$\ be the $(n+1)$-matrix $E_{n}\dotplus(1)$. For an element
$\mathbf{v}\in\mathbf{P}^{n}$, take a column vector $^{t}(s_{1},\cdots
,s_{n+1})$\ of length $1$\ representing $\mathbf{v}$. Let $G(\mathbf{v}%
)$\ denote the $(n+1)$-matrix whose $(i,j)$\ component is $\delta_{ij}%
-2s_{i}s_{j}$, where $\delta_{ij}=1$\ if $i=j$ and $\delta_{ij}=0$\ if $i\neq
j$. Then we set $\mathfrak{g}(\mathbf{v})=I_{-}G(\mathbf{v})$. We note that
$\mathfrak{g}$\ is well-known as the characteristic map of the tangent bundle
of $S^{n+1}$ (see \cite[Section 23.4]{Ste}). We define the space
$\mathcal{M}(n+2)$ to be the union $C_{+}\mathbf{P}^{n}\times SO(n+1)\cup
_{\mathbf{P}^{n}\times SO(n+1)}C_{-}\mathbf{P}^{n}\times SO(n+1)$, which is
pasted by the diffeomorphism $(\mathbf{v},S)\rightarrow(\mathbf{v}%
,\mathfrak{g}(\mathbf{v})S)$ for $S\in SO(n+1)$. We have the map
$\mathfrak{p}_{S\mathbf{P}^{n}}:\mathcal{M}(n+2)\rightarrow S\mathbf{P}^{n}$,
which induces the structure of the fiber bundle with fiber $S(n+1)$.

Let $\mathfrak{p}_{S^{n+1}}:SO(n+2)\rightarrow SO(n+2)/SO(n+1)\times
E_{1}=S^{n+1}$\ be the map defined by $\mathfrak{p}_{S^{n+1}}(S)=S\mathbf{e}%
_{n+2}$\ for $S\in SO(n+2)$.\ Let $D_{+}^{n+1}$\ and $D_{-}^{n+1}$\ denote the
hemispheres of $S^{n+1}$ which consist of\ all points $(x_{1},\cdots,x_{n+2}%
)$\ such that $x_{n+2}\geq0$\ and $x_{n+2}\leq0$ respectively. Let us recall
that $SO(n+2)$ is identified with the union $D_{+}^{n+1}\times SO(n+1)\cup
_{S^{n}\times SO(n+1)}D_{-}^{n+1}\times SO(n+1)$, which is pasted by the
diffeomorphism $(\mathbf{s},S)\rightarrow(\mathbf{s},\mathfrak{g}%
([\mathbf{s]})S)$.\ By the construction we have the natural bundle map
$p_{\mathcal{M}}^{SO}:SO(n+2)\rightarrow\mathcal{M}(n+2)$ covering
$S(p_{\mathbf{P}^{n}})$.%
\begin{equation}%
\begin{array}
[c]{ccc}%
SO(n+2) & \overset{\underrightarrow{\text{ \ \ }p_{\mathcal{M}}^{SO}\text{
\ \ }}}{} & \mathcal{M}(n+2)\\
{\scriptstyle{\mathfrak{p}_{S^{n+1}}}}\downarrow\text{ \ \ \ \ } &  & \text{
\ \ \ }\downarrow{\scriptstyle{\mathfrak{p}_{S\mathbf{P}^{n}}}}\\
S^{n+1} & \overset{\underrightarrow{\text{ \ }S(p_{\mathbf{P}^{n}})\text{ \ }%
}}{} & S\mathbf{P}^{n}%
\end{array}
\end{equation}
According to \cite[Theorem 3.5]{Omega10}\ and \cite[Section 2]{FoldSurg}, we
have the topological embeddings $i_{\mathcal{M}}:\mathcal{M}(n+2)\rightarrow
\Omega^{1}(n+1,n+1)$\ and $i_{n+1}^{SO}:SO(n+2)\rightarrow\Omega
^{(1,0)}(n+1,n+1)$\ such that $i_{\mathcal{M}}(\mathcal{M}(n+2))$\ and
$i_{n+1}^{SO}(SO(n+2))$\ are the deformation retracts of $\Omega^{1}%
(n+1,n+1)$\ and $\Omega^{(1,0)}(n+1,n+1)$\ respectively and that the following
diagram commutes%
\begin{equation}%
\begin{array}
[c]{ccc}%
SO(n+2) & \overset{\underrightarrow{\text{ \ \ \ }p_{\mathcal{M}}^{SO}\text{
\ \ }}}{} & \mathcal{M}(n+2)\\
{\scriptstyle{i_{n+1}^{SO}}}\downarrow\text{ \ \ \ \ } &  & \text{ }%
\downarrow{\scriptstyle{i_{\mathcal{M}}}}\\
\Omega^{(1,0)}(n+1,n+1) & \overset{\underrightarrow{\text{ \ \ }i_{\Omega^{1}%
}^{\Omega^{(1,0)}}\text{ \ \ }}}{} & \Omega^{1}(n+1,n+1),
\end{array}
\end{equation}
where $i_{\Omega^{1}}^{\Omega^{(1,0)}}$ denotes $\pi_{1}^{2}|\Omega
^{(1,0)}(n+1,n+1):\Omega^{(1,0)}(n+1,n+1)\rightarrow\Omega^{1}(n+1,n+1)$.

\begin{lemma}
Let $n+1=4m$. Then $(p_{\mathcal{M}}^{SO})_{\ast}:\pi_{n}(SO(n+2))\rightarrow
\pi_{n}(\mathcal{M}(n+2))$ is injective.
\end{lemma}

\begin{proof}
Consider the following diagram, which is induced from the homomorphism
$(p_{\mathcal{M}}^{SO})_{\ast}$ between the exact sequences for the fiber
bundles $\mathfrak{p}_{S^{n+1}}$ and $\mathfrak{p}_{S\mathbf{P}^{n}}$ in
(7.2).%
\[%
\begin{array}
[c]{ccccc}%
\pi_{n+1}(S^{n+1}) & \overset{\underrightarrow{\text{ \ \ }\partial\text{
\ \ }}}{} & \pi_{n}(SO(n+1)) & \longrightarrow & \pi_{n}(SO(n+2))\\
\downarrow &  & \downarrow &  & \downarrow\\
\pi_{n+1}(S\mathbf{P}^{n}) & \overset{\underrightarrow{\text{ \ \ }%
\partial\text{ \ \ }}}{} & \pi_{n}(SO(n+1)) & \longrightarrow & \pi
_{n}(\mathcal{M}(n+2))
\end{array}
\]
It is known that $\pi_{n}(SO(n+1))\approx\mathbb{Z}\oplus\mathbb{Z}$, $\pi
_{n}(SO(n+2))\approx\mathbb{Z}$\ and $\partial:\pi_{n+1}(S^{n+1}%
)\rightarrow\pi_{n}(SO(n+1))$\ is injective. Then the assertion follows from
Lemma 7.5.
\end{proof}

Consider the homomorphism%
\begin{equation}
(i_{\Omega^{1}}^{\Omega^{(1,0)}})_{\ast}:H^{n+1}(\widehat{V};\pi_{n}%
(\Omega^{(1,0)}(n+1,n+1)))\rightarrow H^{n+1}(\widehat{V};\pi_{n}(\Omega
^{1}(n+1,n+1))),
\end{equation}
which is injective by (7.3) and Lemma 7.6. Therefore, we have the following proposition.

\begin{proposition}
Let $n+1=4m$. Let $V$ be parallelizable. Then, for the obstruction
$\mathfrak{o}(\tau(\widehat{V},\mathcal{T}(f)))$, $(i_{\Omega^{1}}%
^{\Omega^{(1,0)}})_{\ast}(\mathfrak{o}(\tau(\widehat{V},\mathcal{T}(f))))$
becomes the unique obstruction to find a section of $\Omega^{1}(\tau
(\widehat{V},\mathcal{T}(f)),\varepsilon_{\widehat{V}}^{n+2})$.
\end{proposition}

We note that $(i_{\Omega^{1}}^{\Omega^{(1,0)}})_{\ast}(\mathfrak{o}%
(\tau(\widehat{V},\mathcal{T}(f))))$ is not necessarily a primary obstruction.

In the rest of the paper we are only concerned with the case $n<8$. We have
chosen $V$ to be parallelizable. This is justified by the following lemma and
codim$\Sigma^{3}(n+1,n+1)=9$. By the definition of $\tau(\widehat
{V},\mathcal{T}(f))$ in the case $\partial V=S^{n}$, $\mathcal{T}(f)$\ yields
the section of $\Omega^{1}(\tau(\widehat{V},\mathcal{T}(f)),\varepsilon
_{\widehat{V}}^{n+2})|_{CS^{n}}$, which we denote by $s_{\mathcal{T}(f)}$.

\begin{lemma}
Let $n<8$ and $\partial V=S^{n}$. Then we have the following.

(i) If $P_{1}(\tau(\widehat{V},\mathcal{T}(f))$ does not vanish, then any
extension $F:V\rightarrow S^{n}\times\lbrack0,1]$ such that $F|S^{n}\times0=f$
has the singularities of codimension $4$ and of type $(2)$.

(ii) If $P_{1}(\tau(\widehat{V},\mathcal{T}(f))$ vanishes, then we have

(ii-a) $s_{\mathcal{T}(f)}$ is extendable to a section of $\Omega^{1}%
(\tau(\widehat{V},\mathcal{T}(f)),\varepsilon_{\widehat{V}}^{n+2})$ over
$\widehat{V}$\ except for a single point in the interior of $V$. In
particular, $V$ is parallelizable.

(ii-b) $P_{2}(\tau(\widehat{V},\mathcal{T}(f))$ vanishes if and only if
$s_{\mathcal{T}(f)}$ is extendable to a section of $\Omega^{1}(\tau
(\widehat{V},\mathcal{T}(f)),\varepsilon_{\widehat{V}}^{n+2})$ over
$\widehat{V}$.
\end{lemma}

\begin{proof}
(i) is clear.

(ii-a) We triangulate $\widehat{V}$ so that $CV$ is a subcomplex. Since
codim$\Sigma^{2}(n+1,n+1)=4,$ $\Omega^{1}(n+1,n+1)$ is $2$-connected. By
considering the fiber bundle $\mathfrak{p}_{S\mathbf{P}^{n}}:\mathcal{M}%
(n+2)\rightarrow S\mathbf{P}^{n}$ it follows that $TV$ is trivial on $CV$ and
the $3$-skeleton of $V$. Since $P_{1}(\tau(\widehat{V},\mathcal{T}(f))=0$, it
follows from \cite[Proof of Lemma 2]{MiKerBer} that $TV$ is trivial on $CV$
and the $4$-skeleton of $V$. Then the assertion follows from $\pi
_{i}(SO(8))=0$ for $i=4,5,6$ by applying the obstruction theory for
$\mathfrak{p}_{S\mathbf{P}^{n}}:\mathcal{M}(n+2)\rightarrow S\mathbf{P}^{n}$
with fiber $SO(8)$.

(ii-b) We can define $\mathfrak{o}(\tau(\widehat{V},\mathcal{T}(f)))$ by
(ii-a), and hence the assertion follows from Lemma 7.3 and Proposition 7.7.
\end{proof}

Let us see that the singularities of symbol $I([f])$\ of $F$ actually detect
the stable homotopy class $\omega_{0}([f])\in\pi_{n}^{S}$ for $\pi_{1}%
^{S}\approx\pi_{2}^{S}\approx\mathbb{Z}/(2)$, $\pi_{3}^{S}\approx
\mathbb{Z}/(24)$ and $\pi_{7}^{S}\approx\mathbb{Z}/(240)$. We use the above
notation. Note that in dimensions $n=1,2$, stable tangent bundles $\tau_{V}$
and $\tau_{W}$ are trivial (note that an orientable $3$-manifold is parallelizable).

(Case: $n=1$) We may take $N=S^{1}$. We have by \cite[Section 38]{Ste} that
$\mathfrak{o}(\tau_{V},\mathcal{T}(f))=\mathfrak{o}(\tau(\widehat
{V},\mathcal{T}(f)))$ is equal to the second Stiefel\ Whitney class
$w_{2}(\widehat{V})$. This is, as an invariant in $\mathbb{Z}/(2)$, coincides
with the number of the singularities of the symbol $(1,1,0)$ of $F$ modulo
$2$, since the Thom polynomial is the dual class of $S^{\mathbf{1}_{2}}(F)$.
Hence, we have $\mathfrak{m}(1)=0$.

(Case: $n=2$) It follows from Theorem 7.1 that we can choose an $\Omega
^{(1,1,0)}$-regular map $F$ for a fold-map $f$. Hence, $\mathfrak{o}(\tau
_{V},\mathcal{T}(f))$, namely $\mathfrak{o}(\tau(\widehat{V},\mathcal{T}(f)))$
lies in $H^{2}(\widehat{V};\pi_{1}(SO(\ell)))$ and coincides with $w_{2}%
(\tau(\widehat{V},\mathcal{T}(f)))$ by Remark 7.2. Suppose that $w_{2}%
(\tau(\widehat{V},\mathcal{T}(f)))$ vanishes. Since $\pi_{2}(SO(3))\approx
\{0\}$, the second obstruction in $H^{3}(\widehat{V};\pi_{2}(SO(3)))$, for
$\tau(\widehat{V},\mathcal{T}(f))$ to be trivial, always vanishes. This
implies that $\omega_{0}([f])\neq0$\ if and only if $w_{2}(\tau(\widehat
{V},\mathcal{T}(f)))$\ does not vanish for any choice of $V$ and $F$
(consequently, $I([f])=(1,1,0)$).

(Case: $n=3$) By Theorem 7.4 we have that $\mathfrak{m}(3)=j_{1}=24$ and
$a_{1}=2$. By Lemmas 7.3, 7.8 and Proposition 7.7, we have that $\mathfrak{o}%
(\tau(\widehat{V},\mathcal{T}(f)))$ is equal to $\pm P_{1}(\tau(\widehat
{V},\mathcal{T}(f))/2$. By \cite{Ronone} the Thom polynomial of $\Sigma
^{2}(\tau(\widehat{V},\mathcal{T}(f)),\varepsilon_{\widehat{V}}^{4})$ is equal
to $P_{1}(\tau(\widehat{V},\mathcal{T}(f))$. Consequently, an $\Omega^{2}%
$-regular map $F$, for an element $\omega_{0}([f])\in\pi_{3}^{S}%
\approx\mathbb{Z}/(24)$, has the corresponding algebraic number of the
singular points of the symbol $(2,0)$ modulo $24$. We remark that this number
coincides with the $e$-invariant introduced in \cite{JXIV} and \cite{Toda}.

(Case : $n=7$) By \cite[Section 7, Discussions and commputations]{KerMil}, an
element of $\pi_{7}^{S}\approx\mathbb{Z}/(240)$ is detected by $P_{2}%
(\tau(\widehat{V},\mathcal{T}(f))/6$ modulo $240$. We note \textrm{codim}%
$\Sigma^{3}(n+1,n+1)=9$. By Theorem 7.4 we have that $\mathfrak{m}%
(7)=j_{2}=240$ and $a_{1}=1$. Let $f:S^{7}\rightarrow S^{7}$ be a fold map
with $\omega_{0}([f])\neq0$. If $P_{1}(\tau(\widehat{V},\mathcal{T}(f))$ does
not vanish, then we have that $I([f])=(2,0)$ or $(2,1)$ by Lemma 7.8 (i). If
$P_{1}(\tau(\widehat{V},\mathcal{T}(f))$ vanishes, then $V$ becomes
parallelizable by Lemma 7.8 (ii). Since $P_{2}(\tau(\widehat{V},\mathcal{T}%
(f))$\ does not vanish for any cobordism $F:(V,S^{7})\rightarrow(S^{7}%
\times\lbrack0,1],S^{7}\times0)$ with $F|S^{7}=f$, the secondary obstruction
$(i_{\Omega^{1}}^{\Omega^{(1,0)}})_{\ast}(\mathfrak{o}(\tau(\widehat
{V},\mathcal{T}(f))))$ does not vanish by Lemmas 7.3, 7.8 and Proposition 7.7.
Therefore, we have $I([f])=(2,0)$ or $(2,1)$.

Let $IV_{4}=(x^{2}+y^{2},x^{4})$ and $(x^{2}+y^{3},xy^{2})$ stand for the
orbit of the $k$-jets of the $C^{\infty}$-stable germs $(\mathbb{R}%
^{8},0)\rightarrow(\mathbb{R}^{8},0)$ of the symbols $(2,0)$ and $(2,1)$,
which are characterized by the local algebras $\mathbb{R}[[x,y]/(x^{2}%
+y^{2},x^{4})$ and $\mathbb{R}[[x,y]/(x^{2}+y^{3},xy^{2})$, by the group
action of Diff$(\mathbb{R}^{8},0)\times$Diff$(\mathbb{R}^{8},0)$ respectively.
They have been defined in \cite{MathVI}. If we apply an elaborate result in
\cite{FR} to the jet bundle $J^{k}(\tau(\widehat{V},\mathcal{T}%
(f)),\varepsilon_{\widehat{V}}^{8})$, then we obtain the cycle $\left\langle
(x^{2}+y^{3},xy^{2})-2IV_{4}\right\rangle $\ under the integer coefficients of
the Vassiliev complex (\cite[Theorem 2.7]{FR}) and the Thom polynomial of
$\left\langle (x^{2}+y^{3},xy^{2})-2IV_{4}\right\rangle $\ is equal to
$9P_{2}(\tau(\widehat{V},\mathcal{T}(f)))$ (\cite[Section 3]{FR}). We denote
the algebraic numbers of the singular points of types $(x^{2}+y^{3},xy^{2})$
and $IV_{4}$\ by $A$ and $B$ respectively. Then $A-2B$ is divisible by
$6\cdot9=54$ and $(A-2B)/54$ corresponds to the stable homotopy class
$\omega_{0}([f])$.

\end{document}